\pdfoutput=1
\documentclass[reqno]{amsart}

\usepackage[square,numbers,sort&compress]{natbib}
\usepackage{amsmath,amssymb}
\usepackage[foot]{amsaddr}
\usepackage{dblfloatfix}
\usepackage{graphicx}
\usepackage{listings}
\usepackage{xcolor}
\usepackage{tikz,pgfplots}
\usetikzlibrary{intersections}
\usetikzlibrary{arrows}
\usetikzlibrary{decorations.pathmorphing}
\usetikzlibrary{backgrounds}
\usetikzlibrary{fit}
\usetikzlibrary{positioning}
\usetikzlibrary{petri}
\usetikzlibrary{calc}
\usetikzlibrary{through}
\usetikzlibrary{mindmap}
\usetikzlibrary{trees}
\usetikzlibrary{shadows}
\usetikzlibrary{calendar}
\usepackage{tikz-3dplot}
\usetikzlibrary{external}

\usepackage[ruled]{algorithm2e}
\usepackage[utf8]{inputenc}
\newcommand{\fett}[1]{\boldsymbol{#1}}

\makeatletter
\renewcommand{\email}[2][]{%
	\ifx\emails\@empty\relax\else{\g@addto@macro\emails{,\space}}\fi%
	\@ifnotempty{#1}{\g@addto@macro\emails{\textrm{(#1)}\space}}%
	\g@addto@macro\emails{#2}%
}
\makeatother

\title[Streamline method for two-phase flow]{Operator splitting technique using streamline projection for two-phase flow in highly heterogeneous and anisotropic porous media}

\author[E. Vidotto]{Ettore Vidotto$^\dagger{}^\P$}
\author[M. Schneider]{Martin Schneider$^*$}
\author[R. Helmig]{Rainer Helmig$^*$}
\author[B. Wohlmuth]{Barbara Wohlmuth$^\dagger{}^\text{\textdaggerdbl}$}
\address{$^\dagger$ Chair for Numerical Mathematics, Technical University Munich, Germany}
\address{$^\text{\textdaggerdbl}$ Department of Mathematics, University of Bergen, Norway}
\address{$^*$ Department of Hydromechanics and Modelling of Hydrosystems, University of Stuttgart, Germany}
\address{$^\P$ \textnormal{Corresponding author}}
\email{$\lbrace \text{vidotto,wohlmuth}\rbrace$@ma.tum.de}
\email{Martin.Schneider@iws.uni-stuttgart.de, Rainer.Helmig@iws.uni-stuttgart.de}
\date{\today}


\begin{document}

\maketitle

\begin{abstract}
In this paper, we present a fast streamline-based numerical method for the two-phase flow equations in high-rate flooding scenarios for incompressible fluids in heterogeneous and anisotropic porous media. A fractional flow formulation is adopted and a discontinuous Galerkin method (DG) is employed to solve the pressure equation. Capillary effects can be neglected in high-rate flooding scenarios. This allows us to present an improved streamline approach in combination with the one-dimensional front tracking method to solve the transport equation.
To handle the high computational costs of the DG approximation, domain decomposition is applied combined with an algebraic multigrid preconditioner to solve the linear system. Special care at the interior interfaces is required and the streamline tracer has to include a dynamic communication strategy. The method is validated in various two- and three-dimensional tests, where comparisons of the solutions in terms of approximation of flow front propagation with standard fully-implicit finite volume methods are provided.
\end{abstract}

\section{Introduction}
A wide range of applications like nuclear waste storage, drug transport through human tissue or oil recovery involve porous media flow and transport processes. In this paper, we restrict our attention to questions related to the front displacement of the flow in porous media, for which an accurate prediction is important in real-life applications \cite{petruzzelli2013migration, cao2007fractional}. For example, the correct simulation of contaminant leakage into the underground may prevent water sources from being polluted. The simulation of such sharp front phenomena is particularly challenging and standard numerical schemes, like finite difference or finite volume methods, require extra precautions to accurately track the fronts and to handle arbitrary permeability tensors. In this paper, we present an efficient streamline-based method to simulate the flow of a two-phase system of incompressible fluids in highly heterogeneous and anisotropic media, where capillary effects are neglected. If we do not consider the capillary pressure, it corresponds to simulate a scenario where the velocity field does not change drastically with time, therefore causing a negligible diffusion. Such situation may occur if the injection rate is high and is justified for displacements in macroscopic reservoir sections \cite[Chapt. 5.4]{dake2001practice}. Also in \cite{li2004tidal}, capillary effects are neglected for tide-induced groundwater oscillations in coastal aquifers.

The mathematical model that governs the flow of the fluid is provided by the transport equations. The problem consists of a system of partial differential equations (PDE) with pressure, velocity and saturation as unknowns. The governing equations of two-phase flow is written for each phase $\alpha\in\{w,n\}$ as
\begin{equation}\label{eq:saturation_equation}
\frac{\partial (\phi \rho_\alpha S_\alpha)}{\partial t}+\nabla\cdot (\rho_\alpha\fett{v}_\alpha)-\rho_\alpha q_\alpha=0.
\end{equation}
Here, ''$w$'' and ''$n$'' denote a wetting phase and a non-wetting phase, respectively, $\phi$ is the porosity of the medium, $\rho_\alpha$ and $S_\alpha$ are the density and the saturation of phase $\alpha$, $q_\alpha$ represents the source or sink term (e.g., injection or production wells), and $\fett{v}_\alpha$ is the phase velocity defined by the extended Darcy law
\begin{equation}\label{eq:extended_darcy}
\fett{v}_\alpha=-\lambda_\alpha\fett{K}(\nabla p_\alpha-\rho_\alpha \fett{g}),
\end{equation}
where the intrinsic permeability $\fett{K}$ is a symmetric uniformly positive definite tensor. The ratio between relative permeability $k_{r\alpha}=k_{r\alpha}(S_w)$ and dynamic viscosity $\mu_\alpha$ is called phase mobility $\lambda_\alpha=k_{r\alpha}/\mu_\alpha$, $p_\alpha$ is the phase pressure, and $\fett{g}$ is the gravity vector defined as $-g\fett{e}_d$, with gravitational acceleration $g$ and the dimension $d$ of the problem. The two-phase system \eqref{eq:saturation_equation}-\eqref{eq:extended_darcy} can be closed by two additional relations
$$
S_w+S_n=1,\qquad p_c=p_n-p_w,
$$
where the capillary pressure $p_c=p_c(S_w)$ is a function of the wetting phase saturation. As already mentioned, we neglect in this paper the capillary pressure $p_c$, i.e., we set $p_n=p_w$. \\
A mathematically equivalent fractional-flow formulation can be derived (see \cite{helmig1997multiphase}) from the fully-coupled model \eqref{eq:saturation_equation}-\eqref{eq:extended_darcy}. For two incompressible, immiscible fluids and a rigid porous medium, the global pressure fractional flow formulation, in the presence of a source or sink term $q$, can be written as
\begin{equation}\label{eq:fractional_flow_velocity}
\nabla \cdot \fett{v}_t=q, \qquad \fett{v}_t=-\lambda_t(S_w)\fett{K}(\nabla P-\fett{G}),
\end{equation}
\begin{equation}\label{eq:fractional_flow_saturation}
\phi\frac{\partial S_w}{\partial t}+\nabla\cdot (f_w\fett{v}_t)+\nabla\cdot (\lambda_nf_w(\rho_w-\rho_n)\fett{Kg})
=0.
\end{equation}
This system consists of an elliptic pressure equation and a hyperbolic saturation equation. Here, $\fett{v}_t=\sum_\alpha \fett{v}_\alpha$ is the total velocity, which is expressed in terms of the global pressure $P$, the total mobility $\lambda_t=\sum_\alpha\lambda_\alpha>0$, and the term $\fett{G}=\frac{1}{\lambda_t}(\lambda_w\rho_w+\lambda_n\rho_n)\fett{g}$. Furthermore, the fractional-flow function $f_w$ is given by $\lambda_w/\lambda_t$.

The system \eqref{eq:fractional_flow_velocity}-\eqref{eq:fractional_flow_saturation} is solved in this paper using a sequential strategy, where the pressure equation \eqref{eq:fractional_flow_velocity} is first solved to obtain the pressure $P$, and then \eqref{eq:fractional_flow_saturation} is solved using \emph{operator splitting} (OS) techniques, where the advective and the gravitational part are considered separately. A detailed description of the OS idea can be found, e.g., in \cite{holden2010splitting}. OS methods have clear advantages compared to the fully-coupled formulation. For example, the physical and mathematical character of the different terms (advective or diffusive) can be identified. Therefore, a suitable algorithm can be selected for each equation, dependent on the properties and characteristics of the problem. Furthermore, this sequential strategy is very efficient if the pressure and transport equations are weakly coupled, which is satisfied in this paper since capillary effects are neglected. In addition, the unconditionally stable front tracking method is used to solve the advective step, which is a clear advantage compared to the classical \emph{IMPES} (IMplicit Pressure-Explicit Saturation) scheme. On the other hand, fully-coupled implicit methods based, for example, on a discontinuous Galerkin discretization are very stable and also not sensitive to the choice of the time-step size, but may suffer from large numerical diffusion and severe over and undershoot effects \cite{kane2017hp}. In addition, solving the fully-coupled system can be expensive and may require high performance algorithms (see, e.g., \cite{bastian2014fully}).

To solve \eqref{eq:fractional_flow_velocity}, we employ the \emph{symmetric weighted interior penalty Galerkin method} (SWIPG) \cite{di2012analysis, zunino2009discontinuous}. The main advantage of this method consists in its ability to handle discontinuous permeability fields $\fett{K}$, which is a typical property of realistic geological applications, where the coefficient may vary by several orders of magnitude.
Furthermore, the SWIPG method is locally conservative, guarantees high-order accuracy (depending on the regularity of the solution) and can handle full permeability tensors \cite{riviere2008discontinuous, bastian2011benchmark, epshteyn2008convergence}.\\
In particular, an accurate and locally conservative approximation of the velocity field $\fett{v}_t$ is fundamental for the advective step of the sequential algorithm. A combination of the SWIPG numerical solution with an appropriate $H(\mathrm{div})$-projection of the velocity, as in \cite{bastian2003superconvergence}, satisfies all these properties, which are not immediately guaranteed by other standard numerical methods. Algorithms based on finite difference approximations have been used to solve \eqref{eq:fractional_flow_velocity}, but lack of unphysical solutions for highly heterogeneous and anisotropic permeability tensors. Also finite volume methods with a two-point flux approximation (TPFA) yield non-consistent formulations if the permeability tensor $\fett{K}$ is not aligned with the grid directions \cite{aavatsmark2007interpretation}.

Once the velocity field has been obtained, an OS step for solving \eqref{eq:fractional_flow_saturation} is employed, where the advective and the gravity part are treated separately, as presented in \cite{bratvedt1996streamline}. Each part is then solved using a method based on streamline projection \cite{batycky19973d, crane1999streamline, datta2007streamline, vasco1998integrating}. The set of one-dimensional equations along streamlines or gravity lines is then solved by the front tracking method, systematically described, e.g., in \cite{holden2015front, langseth1996implementation, nilsen2009front}.\\
The combination of the streamline method with the front tracking is widely used for simulating subsurface transport \cite{cao2011two, kippe2007method, nilsen2009front}. In fact, streamlines are very efficient to compute and can minimize the numerical diffusion. Furthermore, the front tracking method is very attractive for solving one-dimensional hyperbolic equations due to its unconditional stability, high efficiency and ability to resolve discontinuities. 

In order to increase efficiency for large-scale applications, a domain decomposition on parallel architectures is applied to reduce the computational costs of the SWIPG resolution of \eqref{eq:fractional_flow_velocity}, exploiting the parallelisation possibilities offered by the DUNE framework \cite{bastian2008generic}, where the entire method described in this paper has been developed. Our streamline method is therefore formulated for decomposed domains and requires a dynamic communication strategy between processes. In this paper, a message passing architecture is employed, where each process can only access local memory. A streamline tracer algorithm combined with two-phase flow in porous media on decomposed domain is a novelty. In \cite{bhambri2011compositional}, the authors, proposed an algorithm based on distributed-memory to assign to different processes a part of the streamline set. Thus, this approach does not require any communication between processes for calculating the streamlines. In \cite{camp2011streamline, pugmire2009scalable}, different ways to compute streamlines for visualization purposes on decomposed domains are presented, where the velocity field is already provided from astrophysical or hydraulic simulations. In \cite{gerritsen2009parallel}, the authors proposed a parallel algorithm for single-phase flow on a shared-memory architecture, where all cores have access to the whole memory, and, in particular, to the entire pressure and velocity fields.

This paper is organized as follows: In Section \ref{sect:resolution_method}, we present our method, with particular focus on the communication of our parallel streamline tracer. The method is numerically validated in Section \ref{sect:numerical_experiments}, where two- and three-dimensional results are presented. Conclusions follow in Section \ref{sect:conclusion}.

\section{Numerical Method}\label{sect:resolution_method}
Within this section, a brief description of the numerical methods employed in this paper is provided. In particular, our improved streamline method is presented, which allows us to use larger time steps in \eqref{eq:fractional_flow_saturation}, and still obtain an accurate prediction of the flow front.

\subsection{Solution of the pressure equation}
We illustrate the SWIPG scheme, focusing, in particular, on the choice of the penalty parameter as presented in \cite{bastian2011benchmark,li2015high}. Let us denote with $\Omega\subset\mathbb{R}^d$ the porous medium, where the problem is posed. We subdivide the boundary $\partial\Omega$ into two subsets, $\Gamma_D$ and $\Gamma_N$, where Dirichlet and Neumann boundary conditions are set, respectively. Equation \eqref{eq:fractional_flow_velocity} is completed by the boundary conditions
\begin{equation}\label{eq:def_boundary_conditions}
\begin{aligned}
P&=g_D\qquad &\text{on }\Gamma_D,\\
\fett{v_t}\cdot\fett{n}&=g_N\qquad&\text{on }\Gamma_N.
\end{aligned}
\end{equation}
If the boundary condition is set to be a pure Neumann condition, i.e., $\Gamma_N=\partial\Omega$, then the system is closed by the following compatibility condition for the pressure that guarantees uniqueness
$$
\int_{\partial\Omega}g_N=\int_\Omega q.
$$
Let $\mathcal{E}_h$ be a uniform quadrilateral (in two dimensions) or hexahedral (in three dimensions) mesh of $\Omega$, where $h>0$ is the maximum element diameter and let $\Gamma_h$ denote the set of all interior faces of the mesh. We fix a unit normal vector $\fett{n}_e$ for each interior face $e$ and denote by $E_e^+$ and $E_e^-$ the elements in $\mathcal{E}_h$ such that $e=\partial E_e^+ \cap \partial E_e^-$. With this notation, we set $\fett{n}_e$ to point from $E_e^-$ to $E_e^+$. For a function $v$, we also define its values on both sides of $e$ by $v^\pm_e:=v_{|E_e^\pm}$. The weighted average and jump of a function $v$ on the face $e$ are given by
$$
\{v\}_{e,w}=w_e^- v_e^-+w_e^+ v_e^+,\quad \text{ and } \quad [v]_e=v_e^--v_e^+,
$$
with non-negative weights satisfying $w_e^-+w_e^+=1$. \\
If $e$ is a boundary face, then the average and jump are defined as
$$
\{v\}_{e,w}=v_e^-, \quad \text{ and } \quad [v]_e=v_e^-.
$$
The usual arithmetic average at interfaces corresponds to the particular choice $w_e^+=w_e^-=\frac12$. In this work, we consider a specific permeability-dependent choice for the weights as in \cite{bastian2011benchmark,bastian2014fully,li2015high}. Namely, for all interior faces $e\in\Gamma_h$, we define the weights
$$
w_e^+=\frac{\delta_e^-}{\delta^+_e+\delta_e^-}, \qquad w^-_e=\frac{\delta_e^+}{\delta^+_e+\delta_e^-},
$$
with
$$
\delta_e^\pm=\fett{n}_e^T\cdot\lambda_t(S^\pm)\fett{K}^\pm\cdot\fett{n}_e,
$$
where $S^\pm=S(E_e^\pm)$ are the saturations of the elements $E_e^\pm$.\\
We can now define the SWIPG discretization for the pressure equation, where the DG approximation space is given by 
\begin{equation}\label{eq:dg_basis_space}
V_h^k :=\{v\in L^2(\Omega):v_{|E}\in \mathbb{Q}^k(E)\ \forall E\in\mathcal{E}_h\},
\end{equation}
where $\mathbb{Q}^k=\{p:p=\sum_{\|\alpha\|_\infty\leq k} c_\alpha x^\alpha\}$ in the standard multiindex notation.
In the SWIPG scheme, the discrete solution $P_h\in V^k_h$ satisfies the variational equation
$$
a(P_h,v_h)=\ell(v_h)\qquad\forall v_h\in V_h^k,
$$
with bilinear and linear forms defined, following \cite{oden1998discontinuoushpfinite}, as
$$
\begin{aligned}
a(u,v)=&
 \sum_{E\in\mathcal{E}_h}\int_E \lambda_t\fett{K}\nabla u\cdot\nabla v\\
 &-\sum_{e\in \Gamma_h\cup\Gamma_D}\int_e \{\lambda_t\fett{K}\nabla u \cdot \fett{n}_e\}_{e,w}[v]_e\\
&-\sum_{e\in \Gamma_h\cup\Gamma_D}\int_e \{\lambda_t\fett{K}\nabla v \cdot \fett{n}_e\}_{e,w}[u]_e\\
&+\sum_{e\in \Gamma_h\cup\Gamma_D} \sigma_e \int_e [u]_e[v]_e,\\
\ell(v)=
&\sum_{E\in\mathcal{E}_h}\int_E qv +\lambda_t\fett{K}\fett{G}\nabla v
-\sum_{e\in\Gamma_N}\int_e vg_N\\
&+\sum_{e\in \Gamma_D}\int_e \left( \sigma_e v-\lambda_t\fett{K}\nabla v\cdot\fett{n}_e \right)g_D\\ 
&-\sum_{e\in\Gamma_D\cup \Gamma_h}\int_e \left\{\lambda_t\fett{KGn}_e \right\}_{e,w}[v]_e,\\
\end{aligned}
$$
where $\sigma_e$ is the penalty parameter chosen as in \cite{bastian2011benchmark,bastian2014fully,li2015high}. For each $e\in\Gamma_h$, we define
\begin{equation}\label{eq:penalty_parameter}
\sigma_e=2\beta\frac{\delta_e^+\cdot \delta_e^-}{\delta_e^+ + \delta_e^-}k(k+d-1)\frac{|e|}{\min(|E_e^+|,|E_e^-|)},
\end{equation}
while, for boundary interfaces $e\in \Gamma_D$, we set 
$$
\sigma_e=\beta\delta_ek(k+d-1)\frac{|e|}{|E_e^-|}.
$$
The factor $\beta$ in the penalty term is constant for all faces in our simulations.\\
This method results in a sparse, large, symmetric and positive definite algebraic system of equations for the pressure. These large-scale linear systems (in particular in three dimensions) can be efficiently solved using the parallelisation possibilities offered by the DUNE framework  \cite{bastian2008generic}. Here, we use a conjugate gradient solver together with an AMG preconditioner, see \cite{bastian2012algebraic} for more details.

\subsection{$H(\text{div})$-projection of the velocity}
The velocity field obtained from the DG pressure $P_h$ is not conservative, due to its discontinuities in the normal components of the velocities across element boundaries. To overcome this problem, a post-processing is applied to obtain a conservative velocity field. In this work, we follow the approach presented first in \cite{bastian2003superconvergence} and then employed in other works, e.g., in \cite{lin2015comparative,niessner2007multi}. A projection onto the BDM-space of first order is therefore adopted, which guarantees that the resulting velocity field is continuous across element edges in the normal direction. For quadrilateral and hexahedral elements, the definition of the spaces can be found in \cite{brezzi1985two,brezzi2012mixed}. This projection is an element-wise post-processing and therefore computationally inexpensive.

\subsection{Operator splitting for the transport equation}
For completeness, the OS concept to discretise the saturation equation \eqref{eq:fractional_flow_saturation} is briefly described in this section. For simplicity, let the time interval $I=[0,T]$ be uniformly partitioned into subintervals $I_n = \left(t_{n-1},t_n\right]$ of constant length $\Delta t = t_n-t_{n-1}$. Hereinafter, $S_w^n$ denotes the wetting-phase saturation at time level $n$, and $S_w^0=\mathcal{P}S_{w0}$ the $\mathrm{L}^2$-projection of the initial data $S_{w0}$. By decomposing the spatial differential operator into the advective part $\mathcal{S}_h(\Delta t)$ and the gravitational part $\mathcal{G}_h(\Delta t)$, the OS solution procedure for one splitting step $[t^n,t^{n+1}]$ is defined as
\begin{subequations}
\begin{align}
\phi\frac{\partial S_w}{\partial t}+\nabla\cdot(f_w\fett{v}_t)&=0 & :\ & S_w^n\xrightarrow{\mathcal{S}_h(\Delta t)} \hat{S}_w^n,\label{eq:os_adv_step}\\
\phi\frac{\partial S_w}{\partial t}+\nabla\cdot(\tilde{f}_w\fett{Kg})&=0 & :\ & \hat{S}_w^n\xrightarrow{\mathcal{G}_h(\Delta t)} S_w^{n+1},\label{eq:os_grav_step}
\end{align}
\end{subequations}
with $\tilde{f}_w=\lambda_n f_w(\rho_w-\rho_n)$. $\mathcal{S}_h$ and $\mathcal{G}_h$ represent the discrete solution operators of the advective and the gravitational step, and $\hat{S}_w^n$ is the intermediate saturation value between two calculation steps. Due to the advective character of both differential operators, the front tracking method is applied along streamlines and gravity lines, respectively.

\subsection{Improved streamline method for higher dimensions}
In the following, we present our streamline approach employed to solve equation \eqref{eq:os_adv_step}. The method applies to \eqref{eq:os_grav_step} in the same way. Streamlines are traced along the velocity field using the standard \emph{Pollock method} \cite{pollock1988semianalytical}. This method is provided for orthogonal grids and it assumes that each principal velocity component varies linearly within a cell. The application of this method to the reconstructed velocity field $\fett{v}_t^{\text{BDM}}$ is not straightforward, since this velocity does not vary linearly within an element. Therefore, we consider here a further approximation $\fett{\bar{v}}$ of the velocity field within a cell by computing a weighted average along each edge of the cell:
$$
\fett{\bar{v}}_{|e}=\sum_{i}w_i \fett{v}_t^{\text{BDM}}(\xi_i)\cdot\fett{n}_e,
$$
where $e$ is a face of an element, $\{\xi_i\}$ represents a set of quadrature points in $e$, and $\{w_i\}$ are the corresponding weights such that $\sum_i w_i=1$. Due to the continuity of the BDM-velocity across edges in the normal component, this new velocity $\fett{\bar{v}}$ maintains the same continuity properties. To apply the Pollock method, the velocity $\fett{\bar{v}}$ is then approximated linearly inside the element as described in \cite{pollock1988semianalytical}.\\
A streamline is then described in terms of time-of-flight $\tau$, which represents the time required to travel a distance $s$ along a streamline based on the velocity field $\fett{\bar{v}}$, i.e.,
$$
\tau(s)=\int_0^s\frac{\phi(x)}{|\fett{\bar{v}}(x)|}\;dx.
$$
In order to be able to apply the front tracking method on each streamline, the informations needed to construct the initial function for the Riemann problem have to be collected. For each streamline we record the time necessary to cross other elements, together with the global numeration of those elements and their respective saturation. 
Thus, the full-dimensional transport equation is transformed into a set of one-dimensional equations along the streamlines in terms of time-of-flight. The one-dimensional front tracking can then be applied along each streamline \cite{cao2011two, kippe2007method, nilsen2009front, berre2002streamline}.

When constructing the initial function for the front tracking, the values of the saturation are mapped from the underlying cartesian grid to the streamline grid, introducing some mass balance errors. This problem can be tackled using higher order mapping algorithms \cite{mallison2006improved} or increasing the number of streamlines. In \cite{kippe2007method}, the time-of-flight values are scaled to locally stretch or shrink the streamline grids in order to impose mass-conservation. In this work, we are not interested in forcing mass-conservation, but the focus is set on the accurate approximation of the flow fronts. Therefore, following \cite{bratvedt1996streamline,cao2010robust}, we simply employ a weighted averaging approach for the mapping from the streamline grid to the cartesian one.

\begin{figure}[!ht]
\begin{center}
\begin{minipage}{0.48\textwidth}
\includegraphics[width=1.0\textwidth]{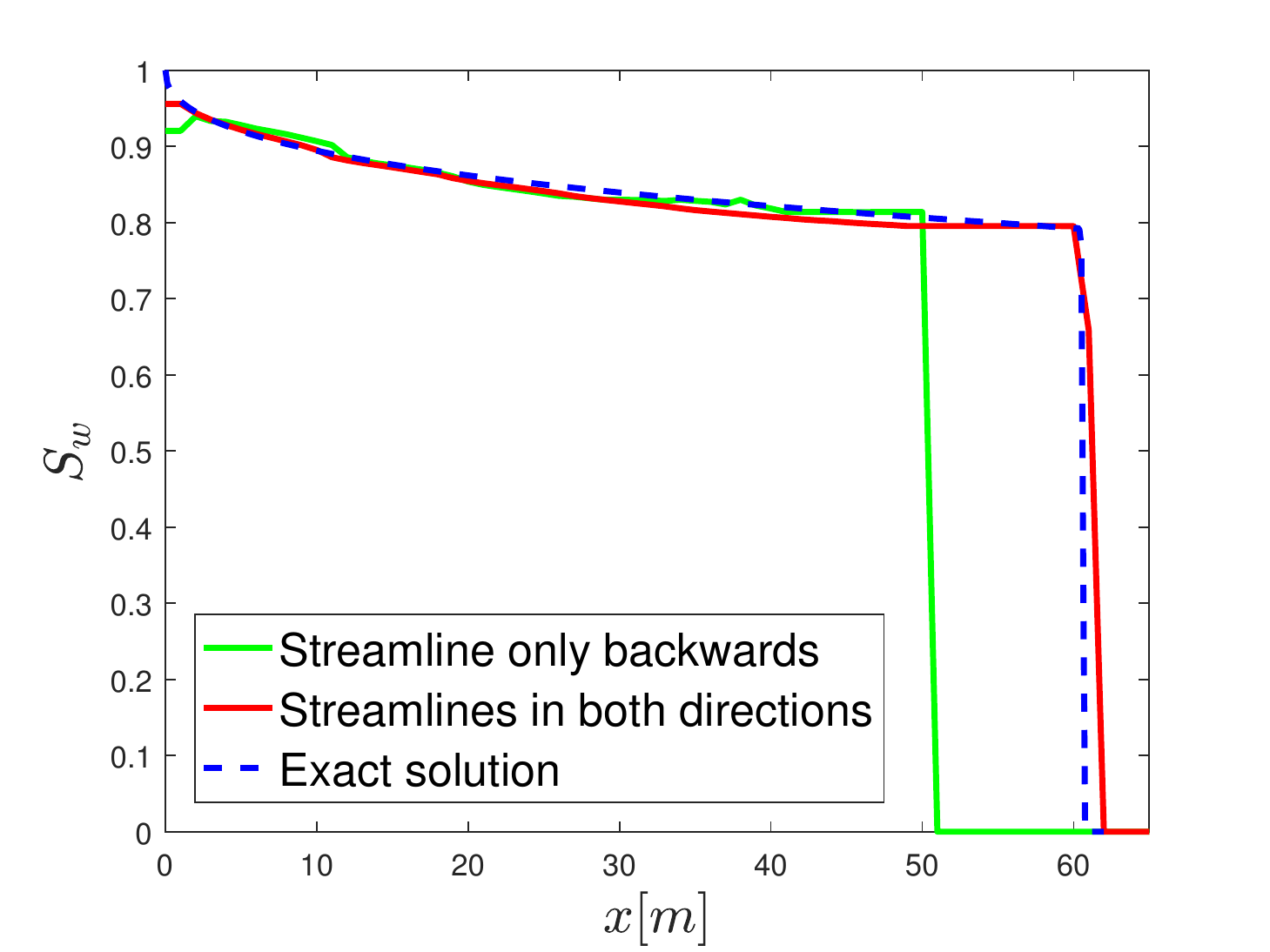}
\end{minipage}
\begin{minipage}{0.48\textwidth}
\includegraphics[width=1.0\textwidth]{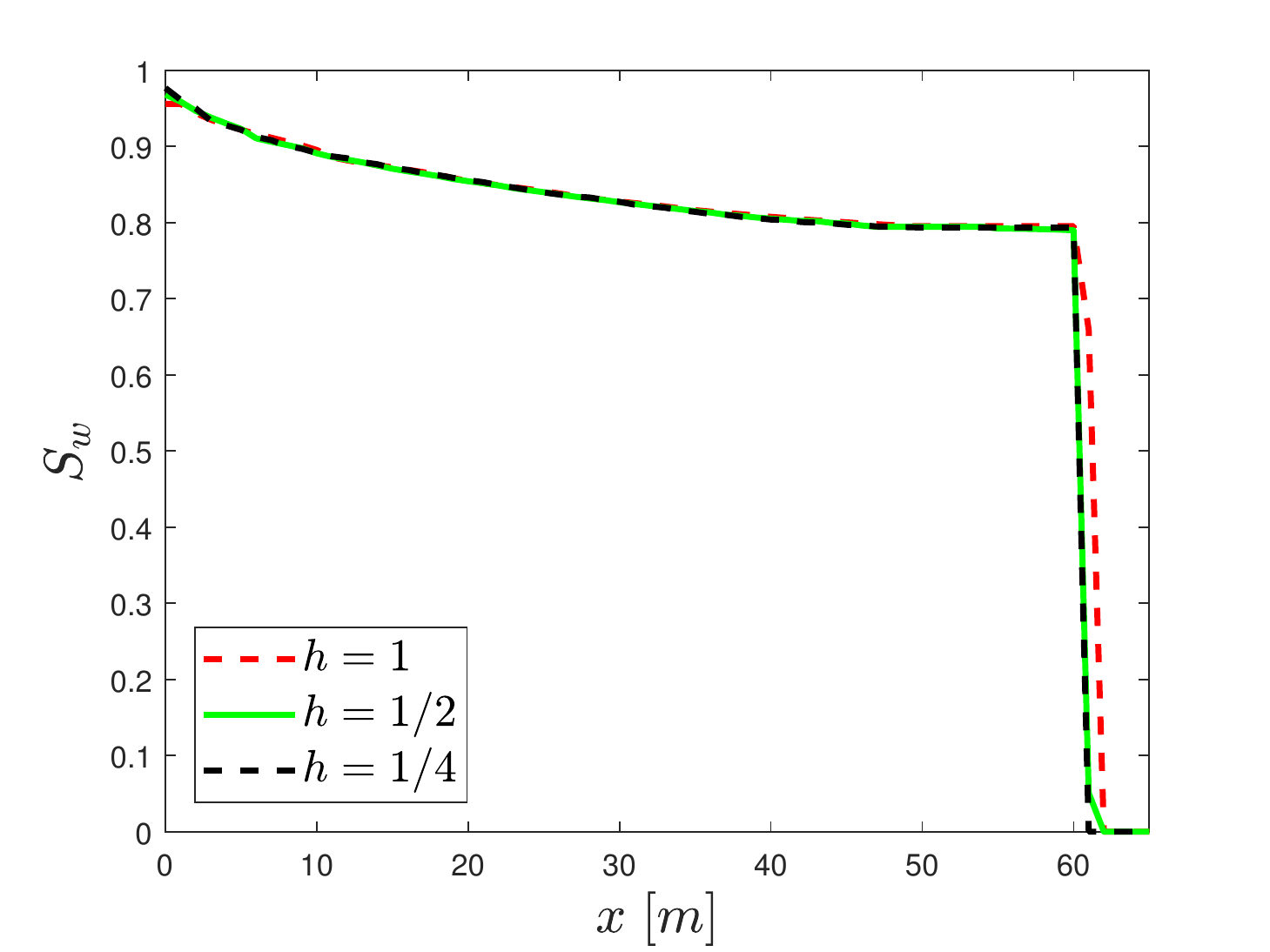}
\end{minipage}
\caption{At the top, the numerical solution of the Buckley-Leverett example is depicted at time $t=10^7\;[s]$. The numerical solution obtained on the mesh with $200\times 200$ elements is compared with the exact solution. At the bottom, the solution fronts are compared with respect to the mesh size \label{fig:comparison_bl}}
\end{center}
\end{figure}

In previous works (see, for example, \cite{ask2000local, borah2013investigation, cao2011two, niessner2007multi, siavashi2014three}), the streamlines are computed backwards along the velocity field for the time interval $\Delta t$, starting from the element centers. Along the backwards part of the streamline, the front tracking method is applied in order to make the transport step forward in time. This approach does not resolve accurately the front propagation, if large time steps are employed (see Fig. \ref{fig:comparison_bl}, top). 
Within this work we therefore apply the front tracking method on the computed streamlines in both directions. This improvement allows us to employ larger time steps and resolve the front propagation more accurately, as shown in the following Buckley-Leverett example.\\
\ \\
\emph{Example.} The difference between the two approaches to solve \eqref{eq:os_adv_step} described above is tested on a simple Buckley-Leverett problem. As computational domain, we consider the square $\Omega=(0,100)^2\subset\mathbb{R}^2$ being initially discretized by a $100\times 100$ mesh. We solve \eqref{eq:fractional_flow_velocity} for the constant permeability $\fett{K}=10^{-10} \fett{I}\;[m^2]$. At the left boundary ($x=0$) we set constant homogeneous Dirichlet boundary conditions $g_D=2\cdot 10^{5}\;[Pa]$ and $S_w=1$, while at the right boundary ($x=100$) Neumann conditions are posed, with $g_N=\frac{1.5\cdot 10^{-3}}{1460}\;[m/s]$. At the top and bottom boundaries, no flow conditions are set, resulting in a constant flux from the left to the right, with zero $y$-component. The mobility function is given by the Brooks-Corey law with $\lambda=2$:
\begin{equation}\label{eq:BrooksCoreyLaw}
\begin{aligned}
k_{rw}(S_w) &= S_w^{\frac{2+3\lambda}{\lambda}}, \\
k_{rn}(S_w) &= (1-S_w)^2\cdot(1-S_w^{\frac{2+\lambda}{\lambda}}).
\end{aligned}
\end{equation}
The viscosities are $\mu_w=10^{-3}\;[Pa\cdot s]$ and $\mu_n=5.7\cdot 10^{-4}\;[Pa\cdot s]$, the porosity is assumed to be constant $\phi=0.2$ for each element. The domain is initially fully saturated with the non-wetting phase, i.e., $S_w=0$ for $t=0$. Gravity effects are neglected in this example. The saturation equation is solved along streamlines using the front tracking method in the time-interval $[0,10^7]$. The difference between both methods can be seen in the top picture of Fig. \ref{fig:comparison_bl}, where the numerical solutions are compared with the analytical one. If only backwards tracing is employed and the time step is large (in this example, the time step has been chosen constant and equal to $\Delta t=10^6\;[s]$), then the propagation front is not accurately resolved. If additionally the streamlines are calculated forwards, then the front matches the exact solution. Furthermore, a relative mass loss of $\approx 16\%$ can be observed, if only backward streamlines are involved, while the total mass loss in the case where streamlines are computed in both directions is reduced to $\approx 0.9\%$. On a uniform refined mesh with $200\times 200$ elements, the relative mass-error is $\approx 0.0661\%$, while on a mesh with $400\times 400$ elements it is $0.0113\%$. We point out that the front is always resolved exactly, as shown in Fig. \ref{fig:comparison_bl}, bottom, independent of the mesh refinement.

\subsection{Parallel implementation of streamlines}
\begin{figure}[!h]
\centering
\begin{tikzpicture}[scale=0.9]
\draw[step=1cm,gray,very thin] (0,0) grid (12,6);
\draw[red,very thick] (0,0) -- (5.98,0) -- (5.98,6) -- (0,6) -- (0,0);
\draw[dotted,blue,very thick] (5,0) -- (5,6);
\draw[blue,very thick] (6.02,0) -- (12,0) -- (12,6) -- (6.02,6) -- (6.02,0);
\draw[dotted,red,very thick] (7,0) -- (7,6);

\node[rotate=90, blue] at (5.5,3) {Overlap P2};
\node[rotate=90, red] at (6.5,3) {Overlap P1};

\node[red] at (3,5.5) {Process 1};
\node[blue] at (9,5.5) {Process 2};
\end{tikzpicture}
\caption{Overlapping domain decomposition using two processes
\label{fig:overlap_process}}
\end{figure}
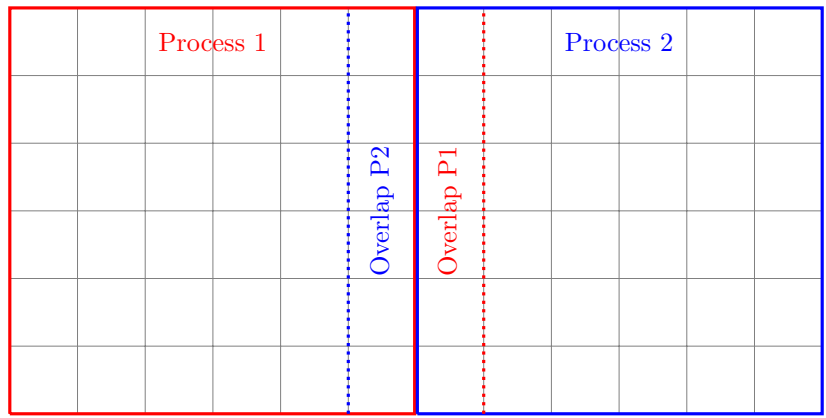

Using the parallel features of DUNE, the pressure field can be computed using an overlapping domain decomposition approach. Therefore, the module for streamline computations has to be parallelized accordingly. A typical situation is depicted in Fig. \ref{fig:overlap_process}, where the computational domain is decomposed into two overlapping subdomains.\\

For each element, streamlines are launched from its centre and are distributed among different processes. Each process starts the computation of its own local set of streamlines independently. If a streamline reaches the boundary of the process where it started, its computation is stopped. Once each process is done with its own set of streamlines, a communication between processes is required to continue the streamlines that have been interrupted at the process boundary. In the following, we present the algorithm employed to track the streamline between different processes. The communication is achieved using \emph{Message Passing Interface} (MPI). Considering the situation presented in Fig. \ref{fig:overlap_process}, each process contains a subdomain, which overlaps with the other process. Between an overlap element in one process and the corresponding interior element in the other process, data can be communicated using the DUNE class:
\begin{lstlisting}
Dune::CommDataHandleIF< DataHandleImp, DataTypeImp >
\end{lstlisting}
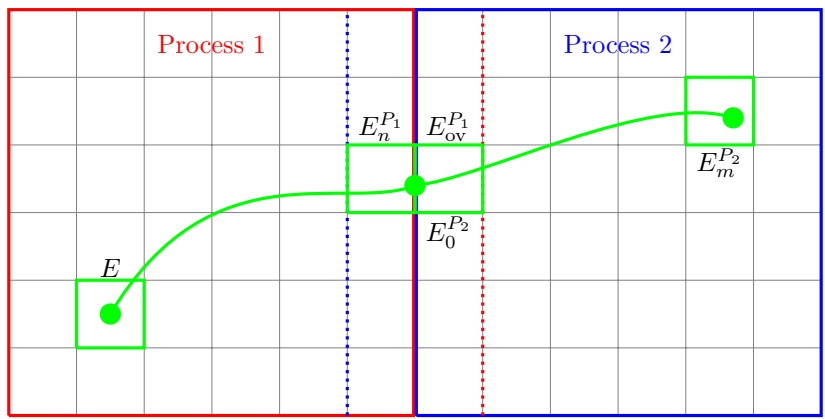
\begin{figure}[!h]
\centering
\begin{tikzpicture}[scale=0.9]

\draw[step=1cm,gray,very thin] (0,0) grid (12,6);
\draw[red,very thick] (0,0) -- (5.98,0) -- (5.98,6) -- (0,6) -- (0,0);
\draw[dotted,blue,very thick] (5,0) -- (5,6);
\draw[blue,very thick] (6.02,0) -- (12,0) -- (12,6) -- (6.02,6) -- (6.02,0);
\draw[dotted,red,very thick] (7,0) -- (7,6);

\node[red] at (3,5.5) {Process 1};
\node[blue] at (9,5.5) {Process 2};

\draw[green,very thick] (1,1) -- (2,1) -- (2,2) -- (1,2) -- (1,1);
\node[shape=circle,inner sep=1mm,fill=green,opacity=0.0,label={[black]above:$E$ }] at (1.5,1.75) {};
\node[shape=circle,inner sep=1mm,fill=green,opacity=1.0] at (1.5,1.5) {};

\draw[green,very thick] (5,3) -- (6,3) -- (6,4) -- (5,4) -- (5,3);
\node[shape=circle,inner sep=1mm,fill=green,opacity=0.0,label={[black]above:$E_n^{P_1}$ }] at (5.5,3.75) {};

\draw[green, very thick] (1.5,1.5) .. controls (3,4) and (5,3) .. (6,3.4);

\draw[green,very thick] (6,3) -- (7,3) -- (7,4) -- (6,4) -- (6,3);
\node[shape=circle,inner sep=1mm,fill=green,opacity=0.0,label={[black]above:$E_{\mathrm{ov}}^{P_1}$ }] at (6.5,3.75) {};
\node[shape=circle,inner sep=1mm,fill=green,opacity=0.0,label={[black]below:$E_0^{P_2}$ }] at (6.5,3.25) {};
\node[shape=circle,inner sep=1mm,fill=green,opacity=1.0] at (6,3.4) {};

\draw[green,very thick] (10,4) -- (11,4) -- (11,5) -- (10,5) -- (10,4);
\node[shape=circle,inner sep=1mm,fill=green,opacity=0.0,label={[black]below:$E_m^{P_2}$ }] at (10.5,4.25) {};

\draw[green, very thick] (6,3.4) .. controls (7,3.5) and (9.5,4.8) .. (10.7,4.4);
\node[shape=circle,inner sep=1mm,fill=green,opacity=1.0] at (10.7,4.4) {};

\end{tikzpicture}
\caption{Streamline crossing the process border 
\label{fig:parallelStreamlineExample} }
\end{figure}
For ease of presentation, let us consider again the simplified problem involving only two processes. We start a streamline from an element $E\in\mathcal{E}_h$ in the first process, as depicted in Fig. \ref{fig:parallelStreamlineExample}. Let us denote by $E_n^{P_2}$ the element where the streamline ends after $\Delta t$ seconds. To be determined are the elements crossed by the streamline, the corresponding crossing times and, eventually, the different saturation values. 
Let the velocity field $\fett{\bar{v}}$ be given. In Algorithm \ref{alg:overlapping_streamline}, a simplified version of the algorithm is presented, where the procedure for tracking the streamlines over a time $\Delta t$ is provided for the case depicted in Fig. \ref{fig:parallelStreamlineExample}. The extension to more processes follows the same concept. Furthermore, we point out that the presented algorithm is independent of the dimension $d$ of the problem.
\begin{algorithm}[!h]
\KwData{$\Delta t,\;\fett{\bar{v}}$}
\KwResult{Streamline of length $\Delta t$ in terms of time of flight;\\
$\qquad\qquad$Elements crossed;\\
$\qquad\qquad$Values of saturation along the streamline.}

\For{Each interior element $E$ on each process}
{
$\bullet$ Compute the streamline from the center of $E$ until time $\Delta t$ or boundary $\partial\Omega$ is reached\;
$\bullet$\If{Streamline reaches an overlap element $E^{P_1}_{\mathrm{ov}}$}
		     {
		     $\bullet$ Interrupt the streamline\;
		     $\bullet$ Save current informations in $E^{P_1}_{\mathrm{ov}}$\;
		     }
}
$\bullet$ Communicate interrupted streamlines from every overlap element $E^{P_1}_{\mathrm{ov}}$ to the corresponding interior element $E^{P_2}_1$\;
\For{Each interrupted streamline starting from an element $E^{P_2}_1$}
{
$\bullet$ Continue the streamline on process $P_2$ until time $\Delta t$ or boundary $\partial\Omega$ is reached\;
$\bullet$ Save informations in $E^{P_2}_1$\;
}
$\bullet$ Communicate informations back from $E^{P_2}_1$ to $E^{P_1}_{\mathrm{ov}}$\;
\For{Each interrupted streamline started in $P_1$}
{
$\bullet$ Return informations from $E^{P_1}_{\mathrm{ov}}$ to the original element $E$.
}
\caption{Pseudo-code for parallel computation of streamlines on decomposed domains \label{alg:overlapping_streamline}}
\end{algorithm}

\section{Numerical experiments}\label{sect:numerical_experiments}
In this section, we validate our method on well-known two- and three-dimensional problems
. We introduce a new parameter, which will be used later in some numerical experiments to validate our method. In the following problems, we assume that the wetting phase is injected from a part $\Omega_I$ of the domain $\bar{\Omega}$, and extracted from another part $\Omega_E$. We assume that $\Omega_I$ is initially filled by the wetting phase, i.e., $S_w(x)=1$ for every $x\in\Omega_I$. On the other end, we assume $S_w(x)=0$ for every $x\in\Omega_E$. For the discretization of problem \eqref{eq:fractional_flow_velocity}-\eqref{eq:fractional_flow_saturation}, we introduce a uniform partition of the time-interval $[0,T]$ into subintervals of length $\Delta t$, where $m\Delta t=T$. We define the \emph{detection time} as the quantity $T_{d}=k\Delta t$, where $k\in\{1,...,m\}$, such that there is at least a $x\in \Omega_E$ with $S_w(x)>0$ and for each $\tilde{k}<k$ the wetting phase has still not reached the extractor, i.e., $S_w(x)=0$ for every $x\in\Omega_E$ at time $\tilde{k}\Delta t$. The actual arrival time of the wetting phase to $\Omega_E$ lies therefore in the interval $(T_d-\Delta t, T_d]$. \\
For all examples, $\mathbb{Q}^1$-elements have been chosen for solving the pressure equation \eqref{eq:fractional_flow_velocity}, i.e., $k=1$ in \eqref{eq:dg_basis_space}, and the parameter $\beta$ in \eqref{eq:penalty_parameter} is set to be constant and equal to one.
\subsection{Five-Spot problem}
The setting for the following problem is the same as in \cite{schneider2016monotone}. As simulation domain, the square $\Omega=(0,100)^2$ is chosen, where the boundary is subdivided in the following subsets:
$$
\begin{aligned}
\Gamma_D&=\{(0,y):y\in[95,100]\}\cup \{(x,100):x\in[0,5]\};\\
\Gamma_N&=\{(100,y):y\in[0,5]\}\cup \{(x,0):x\in[95,100]\};\\
\Gamma^{\mathrm{nf}}_N&=\partial\Omega\setminus(\Gamma_D\cup\Gamma_N).\\
\end{aligned}
$$
On $\Gamma_D$, Dirichlet boundary conditions $g_D=2\cdot 10^5\;[Pa]$ and $S_w=1$ are set, while on $\Gamma_N$ the total velocity in normal direction is given as $g_N=\frac{1.5\cdot 10^{-3}}{1460}\;[m/s]$. On $\Gamma^{\mathrm{nf}}_N$, no-flow condition $g_N=0$ is imposed. Gravity is neglected in this example and the viscosities are given by $\mu_w=10^{-3}\;[Pa\cdot s]$ and $\mu_n=5.7\cdot 10^{-4}\;[Pa\cdot s]$. The porosity is chosen constant in the entire domain, $\phi=0.2$. The relative permeabilities are chosen accordingly to the Brooks-Corey law \eqref{eq:BrooksCoreyLaw} with $\lambda=2$, while the permeability is chosen to be constant $\fett{K}=10^{-10} \fett{I}\;[m^2]$. The simulation interval is $[0,8\cdot 10^7]$, with time step size $\Delta t=5\cdot 10^6\;[s]$. The domain is first discretized by $100\times 100$ elements, yielding a spacial size of $dx=dy=1\;[m]$. A finer mesh is also considered, where the coarse one is uniformly refined, resulting in $dx=dy=1/2\;[m]$.\\
\begin{figure*}[!t]
\begin{center}
\includegraphics[width=0.31\textwidth]{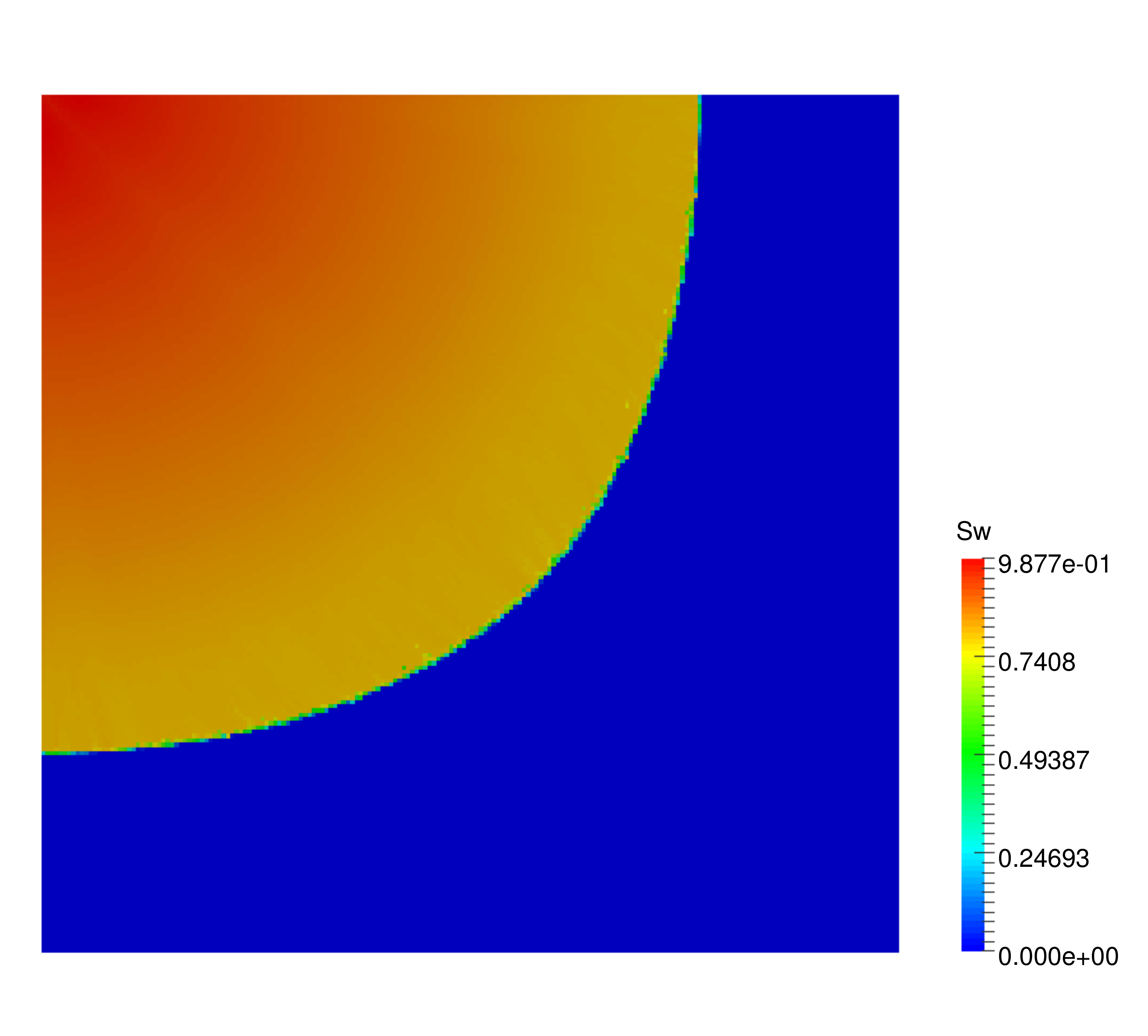}
\includegraphics[width=0.32\textwidth]{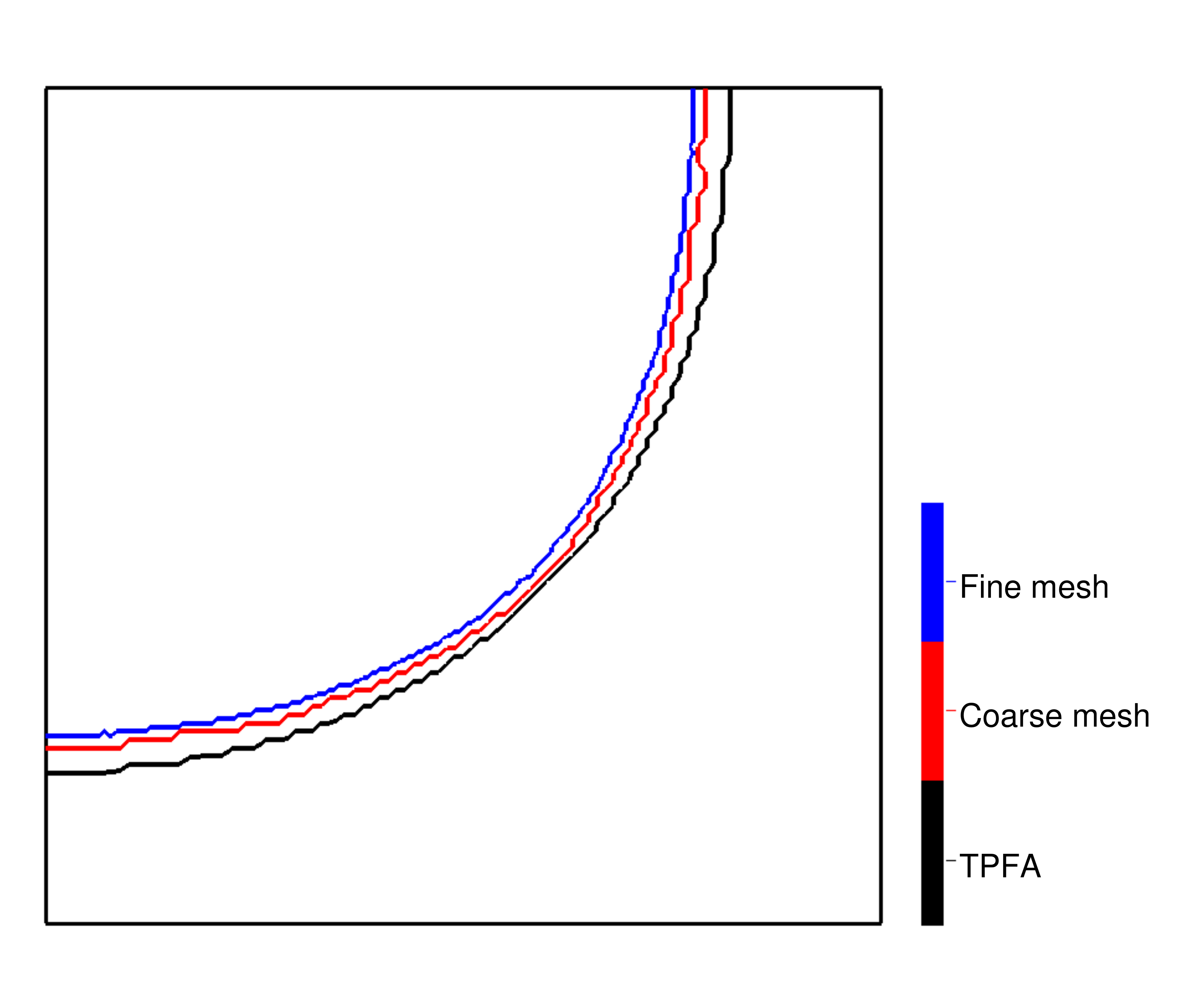}
\hspace{-3mm}
\includegraphics[width=0.36\textwidth]{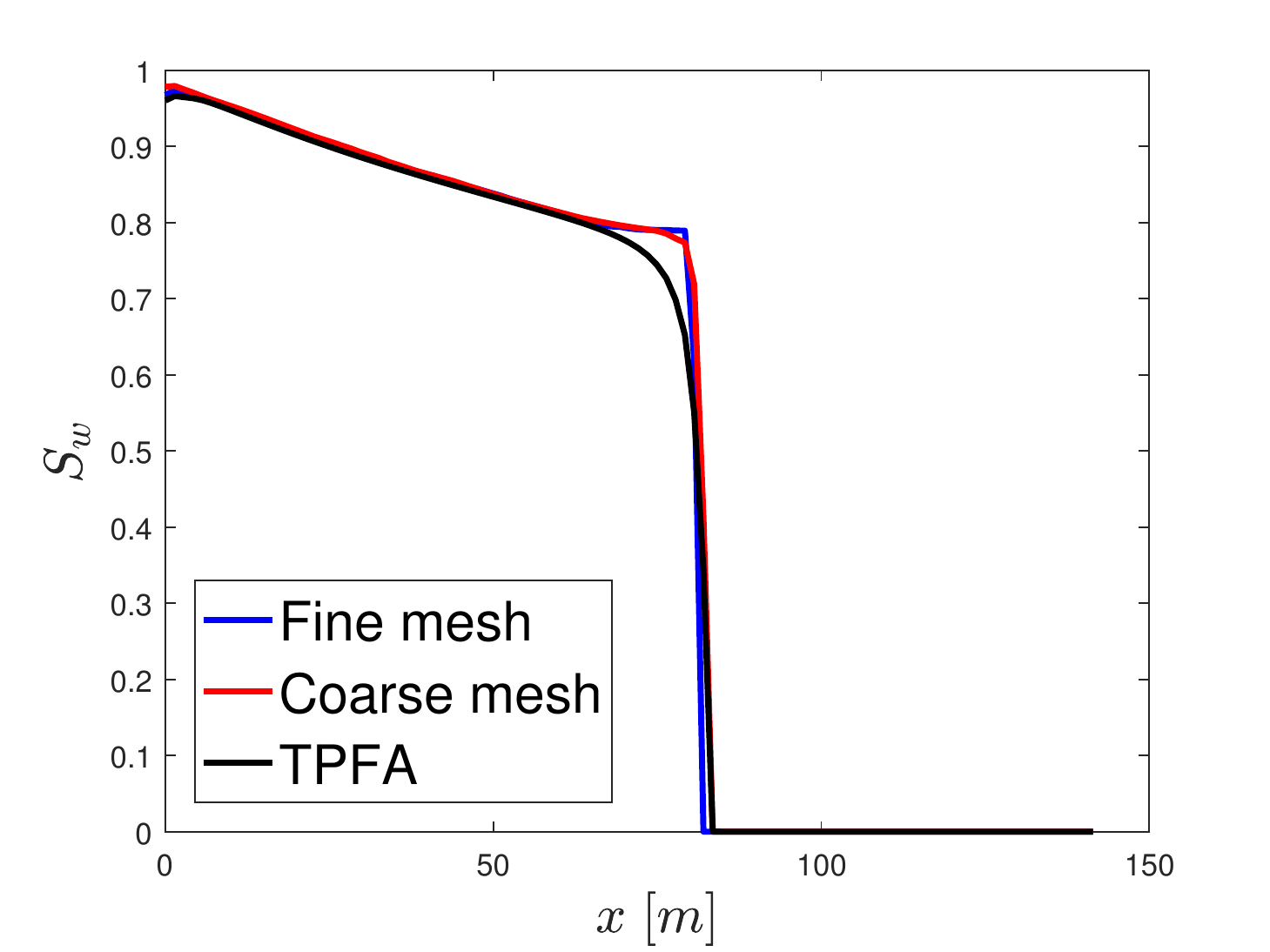}
\caption{The solution obtained with our method at the end of the simulation on the fine grid is shown in the left picture. In the middle, a contour plot for $S_w=10^{-5}$ at end time is provided for the solutions on the coarse and on the fine grid obtained with our method, together with a reference solution computed with a standard fully-implicit TPFA method. The right figure depicts the saturation profile along the diagonal $(0,100)-(100,0)$ \label{fig:5spot_solutions}}
\end{center}
\end{figure*}
The numerical solution obtained with our method at the final time is depicted in Fig. \ref{fig:5spot_solutions}, on the left. In the middle of Fig. \ref{fig:5spot_solutions}, a contour plot for $S_w=10^{-5}$ is shown, where the solutions of our method on both meshes are compared with a reference solution computed with a standard fully-implicit TPFA method on the fine mesh. The reduced numerical diffusion of our method can be observed. \begin{figure*}[!b]
\begin{center}
\includegraphics[width=0.28\textwidth]{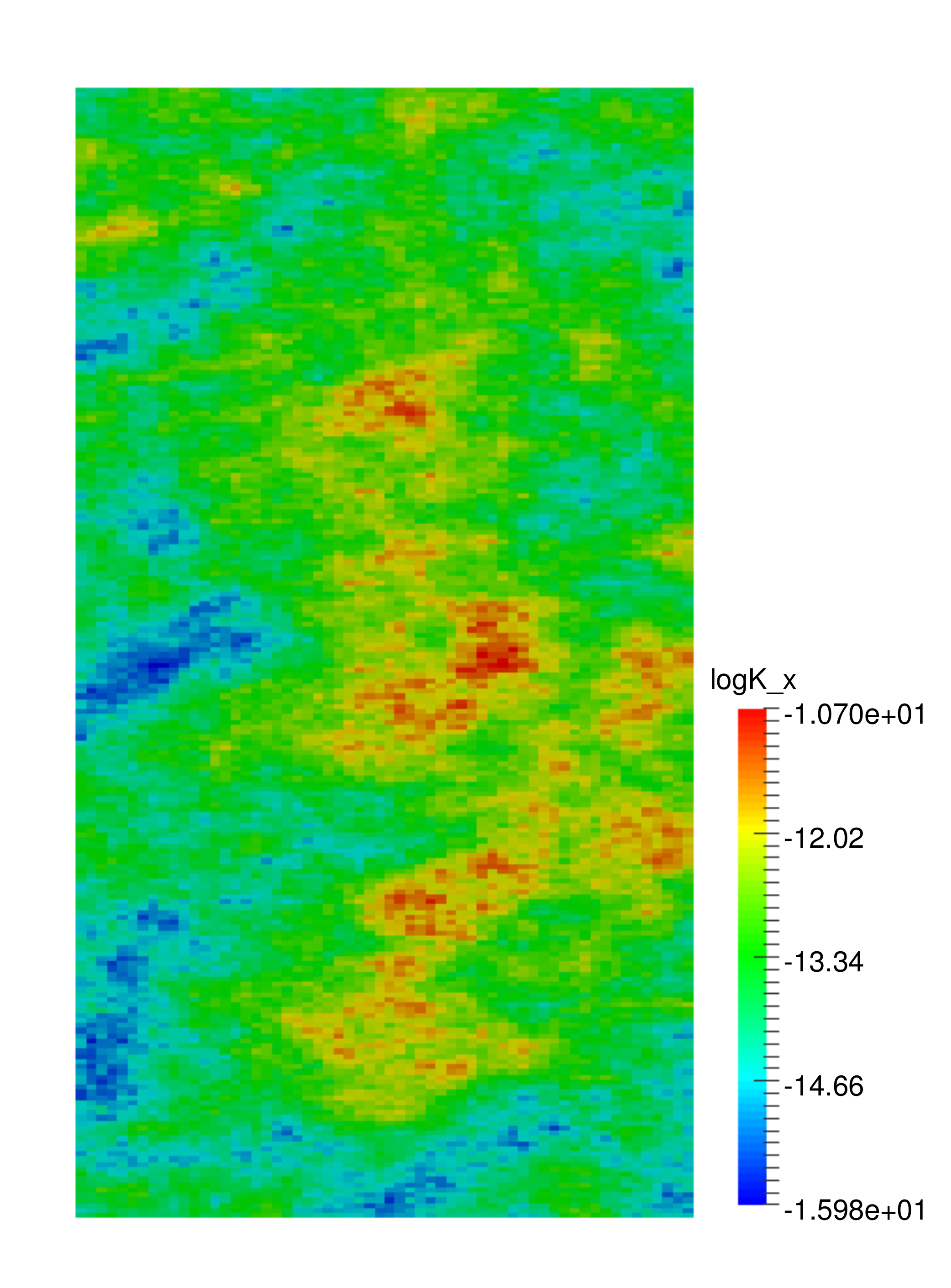}
\includegraphics[width=0.28\textwidth]{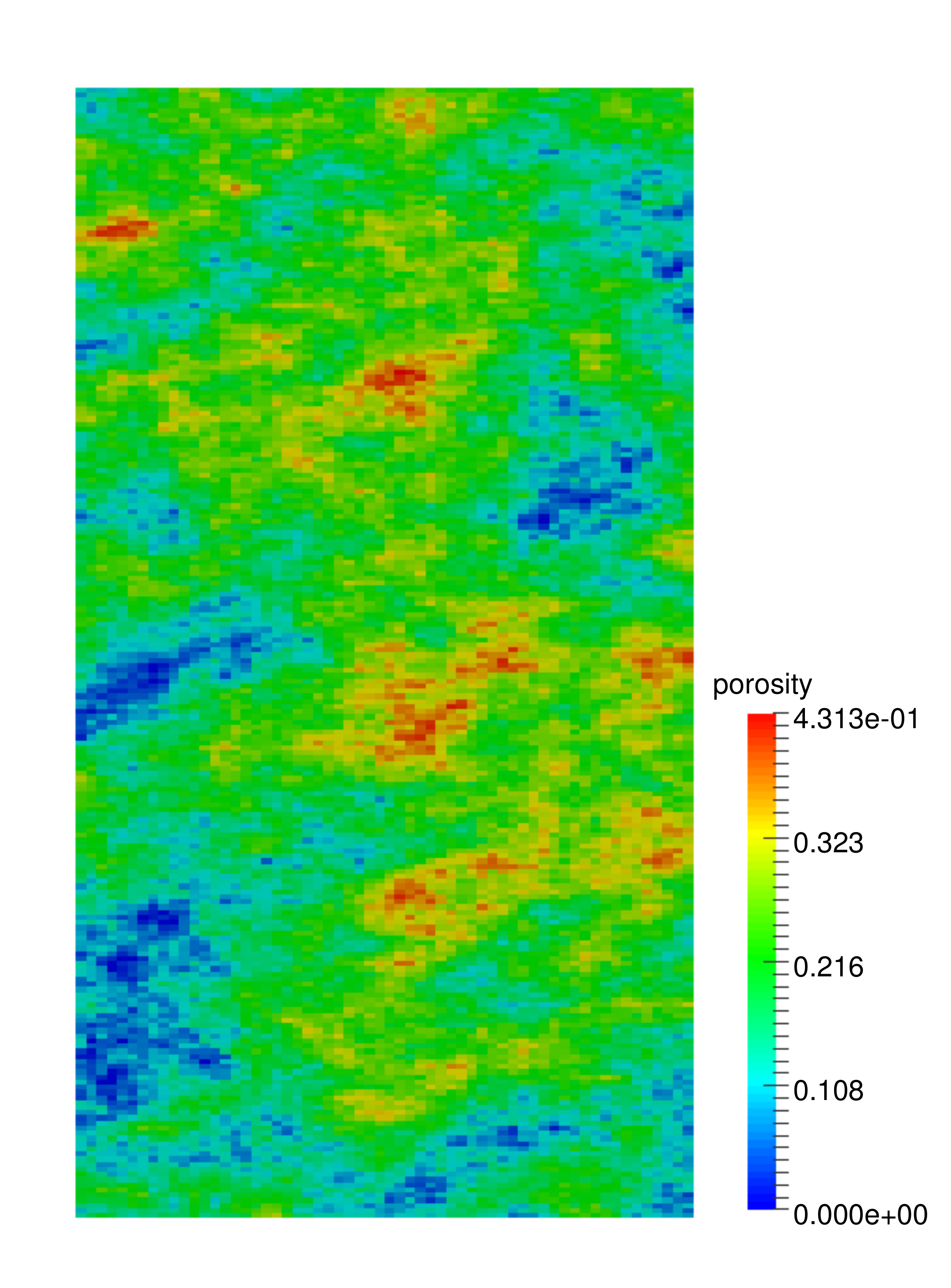}
\includegraphics[width=0.32\textwidth]{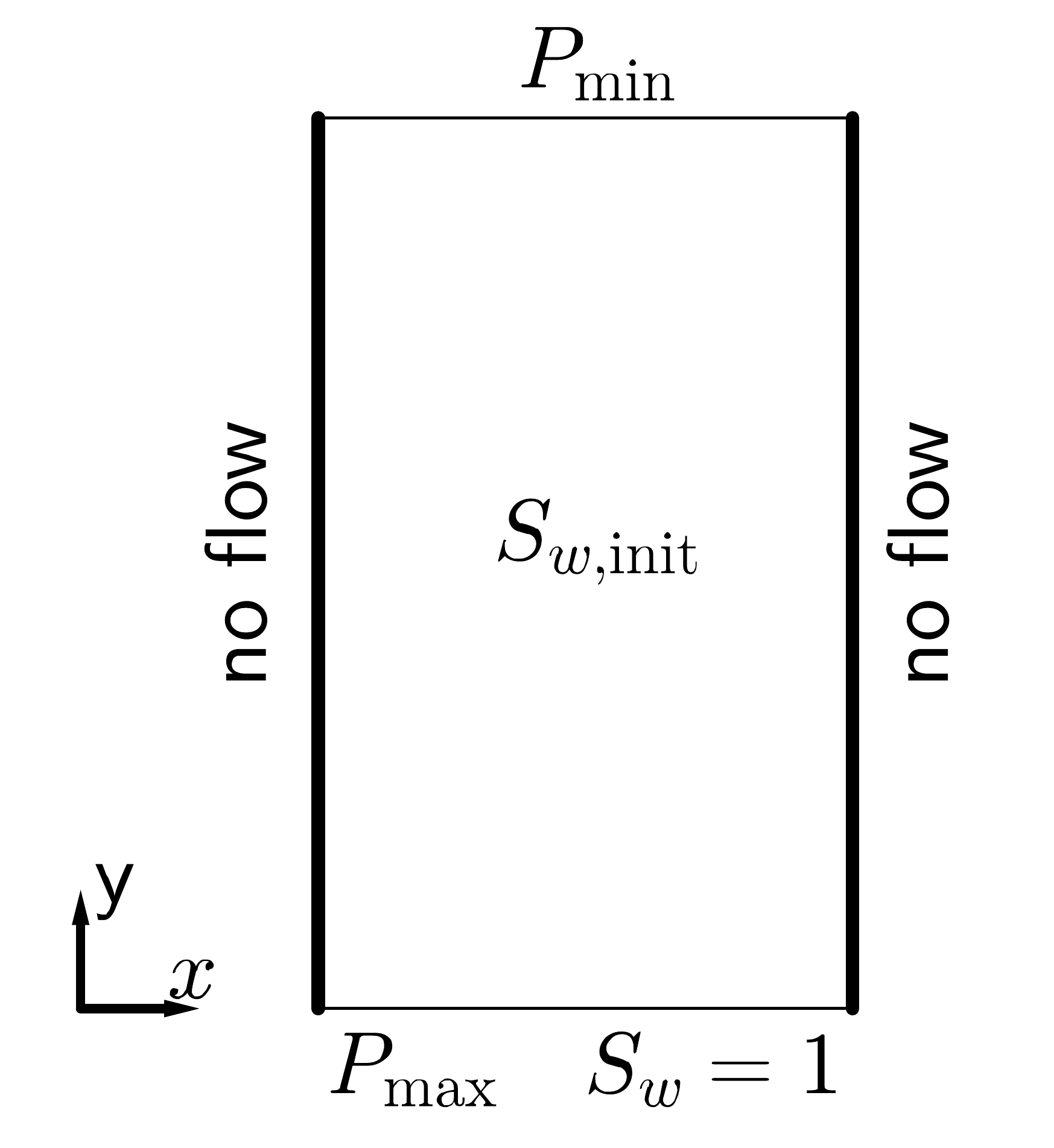}
\caption{In the left picture, the first component of the permeability field is depicted (the medium is isotropic). In the middle picture, the porosity field is shown. On the right, the problem setup is provided\label{fig:spe10_problem_setup}}
\end{center}
\end{figure*}
Furthermore, the wetting phase saturation is plotted in Fig. \ref{fig:5spot_solutions}, on the right, along the diagonal together with the fully-implicit TPFA solution. For the solutions computed with our method, we can observe that numerical diffusion is reduced on the finer mesh and that the two fronts coincide, showing that the front is well resolved by our method, independently of the mesh refinement.

\subsection{Two-dimensional heterogeneous problem}\label{sect:2d_heterogeneous_prob}
In the second example, we test our method on the layer 16 (top formation) of the SPE10 benchmark study \cite{christie2001tenth}. The permeability and porosity fields are shown in Fig. \ref{fig:spe10_problem_setup} in the left and middle pictures, respectively. Both fields show high parameter contrasts. The relative permeabilities are calculated using quadratic laws
\begin{figure*}[!t]
\begin{center}
\includegraphics[width=0.32\textwidth]{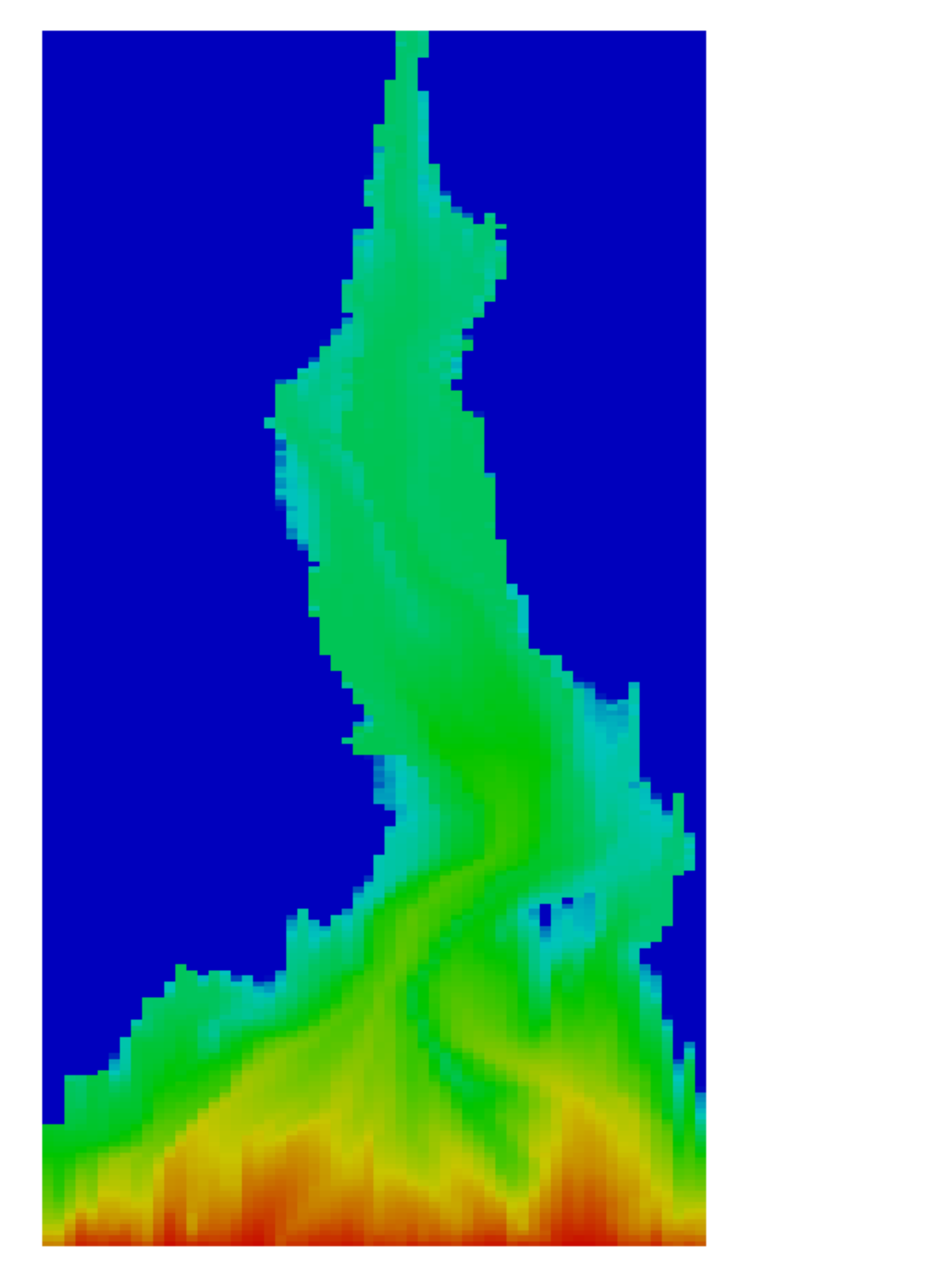}
\includegraphics[width=0.32\textwidth]{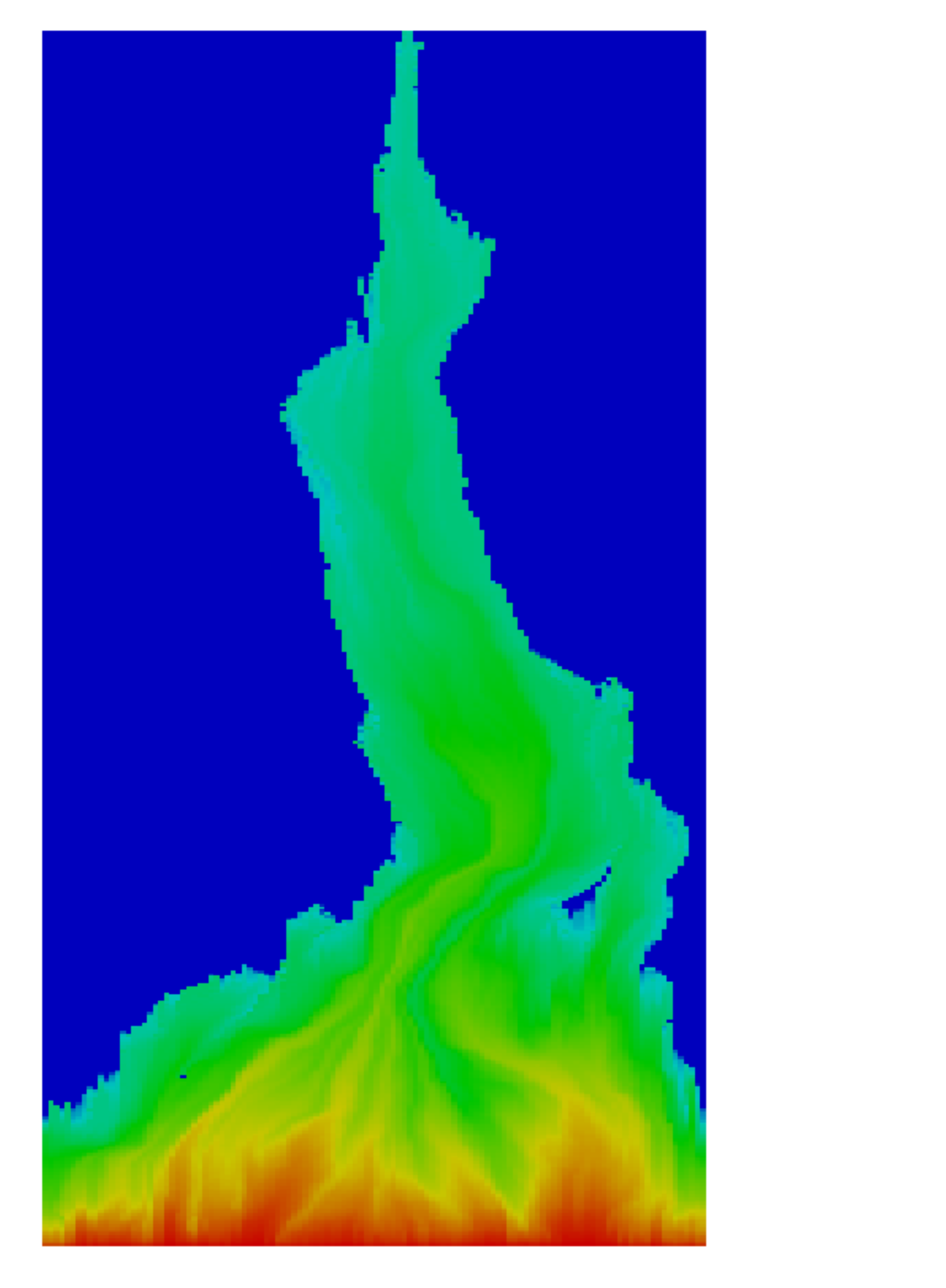}
\includegraphics[width=0.32\textwidth]{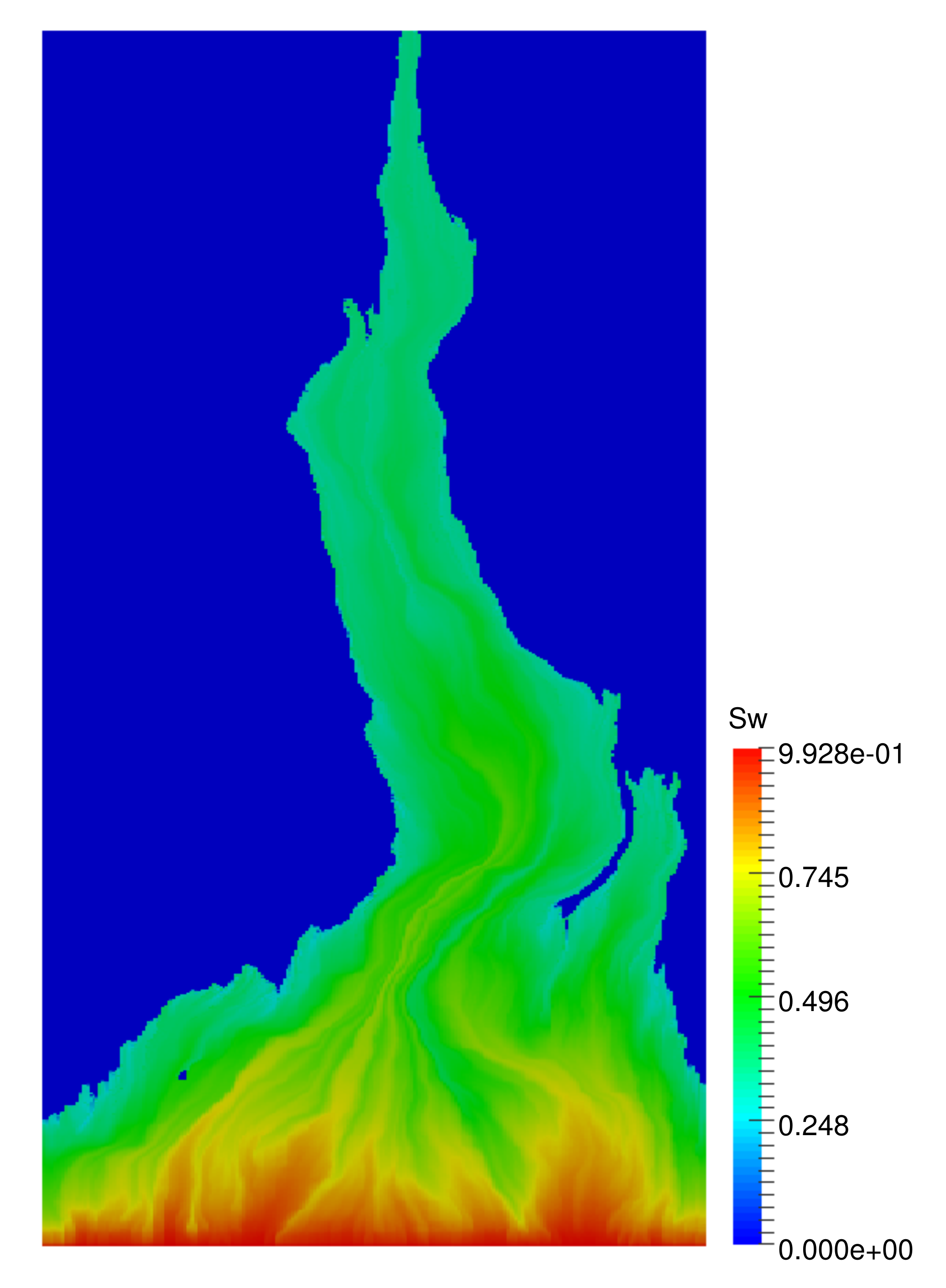}
\end{center}
\end{figure*}
\begin{figure*}[!t]
\begin{center}
\includegraphics[width=0.32\textwidth]{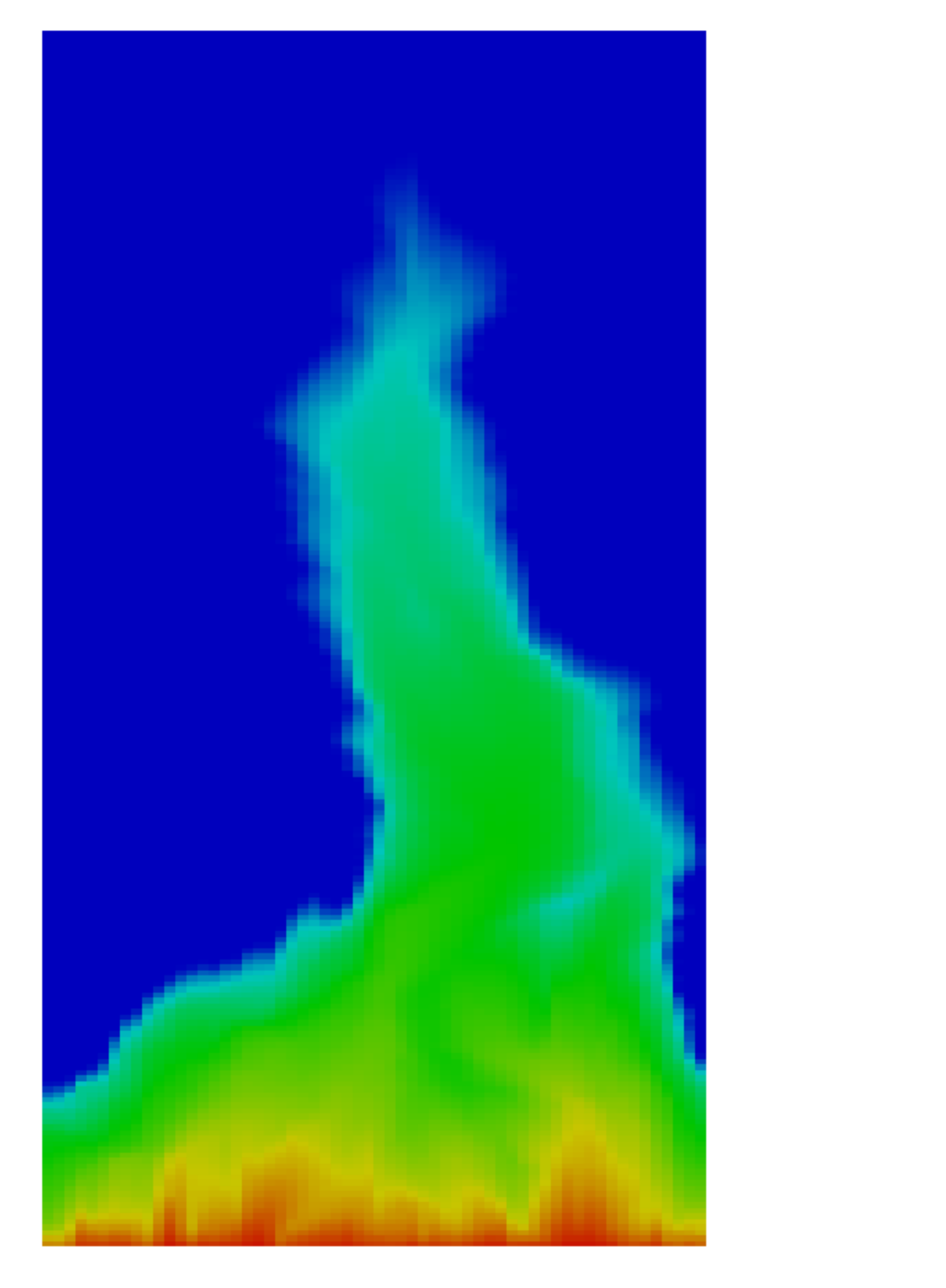}
\includegraphics[width=0.32\textwidth]{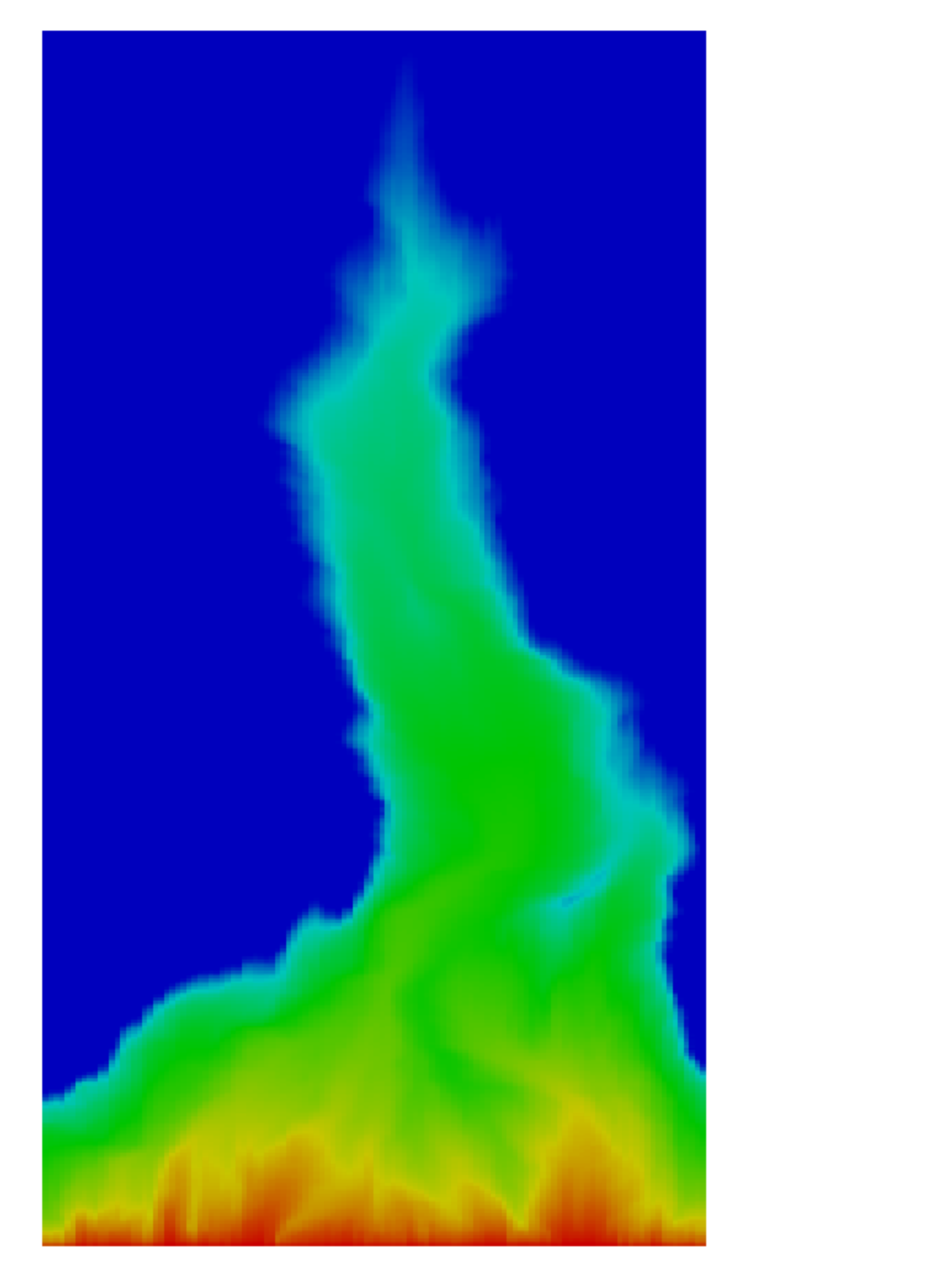}
\includegraphics[width=0.32\textwidth]{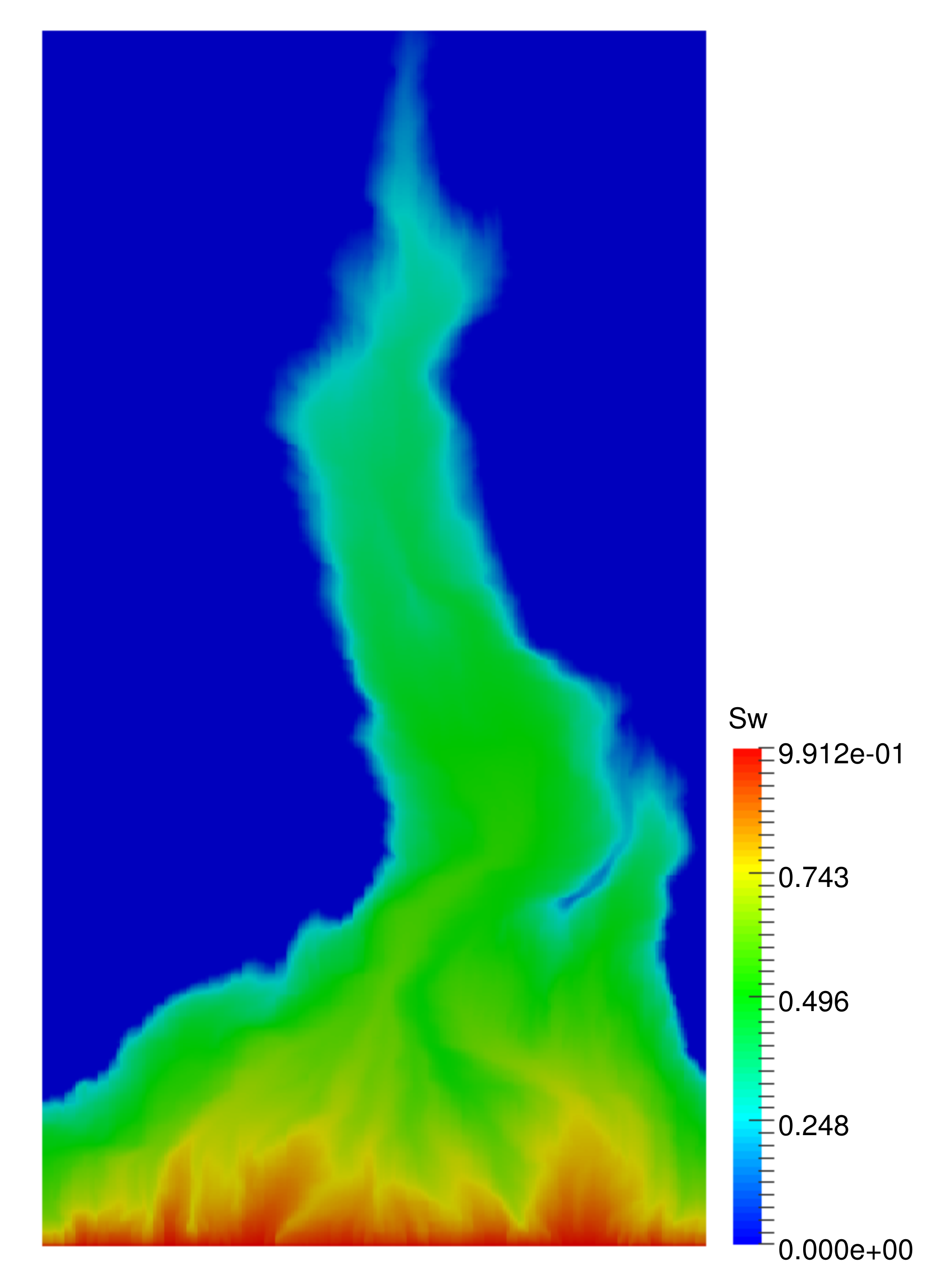}
\caption{At the top, the solutions after $1.3\cdot 10^8\;[s]$ computed with our method on different meshes, from the coarsest (left picture) to the finest (right picture). At the bottom, the reference solutions computed with a fully-implicit TPFA method on the corresponding meshes\label{fig:spe_solutions}}
\end{center}
\end{figure*}
\begin{equation}\label{eq:quadratic_rel_perm}
k_{rw}= S_w^2,\qquad k_{rn}= (1-S_w)^2,
\end{equation}
and the fluid viscosities are $\mu_w=10^{-3}\;[Pa\cdot s]$ and $\mu_n=5\cdot 10^{-3}\;\;[Pa\cdot s]$.
The problem setup is shown in Fig. \ref{fig:spe10_problem_setup}, on the right. The domain is initially saturated by the non-wetting phase (oil), and the wetting phase (water) infiltrates the domain from the lower boundary, i.e., $S_w=1$. The left and right boundaries are closed and a pressure difference of $2\cdot 10^7\;[Pa]$ between the lower and upper boundary is applied. The domain is initially discretized by a grid of $60\times 220$ cells of sizes $dx=6.096\;[m]$ and $dy=3.048\;[m]$. Gravity is neglected. A time step of $\Delta t=10^7\;[s]$ is used for the simulations.
The results are shown in Fig. \ref{fig:spe_solutions}, where the saturation $S_w$ is shown at time $1.3\cdot 10^8\;[s]$. The three solutions at the top have been computed with the method developed in this paper, while the three at the bottom with a standard TPFA method and are referred to as reference solutions for this example. The solutions in the left column have been computed on the initial mesh, while the mesh in the middle column has been uniformly refined and contains $120\times 440$ elements. In the right column, a further uniform refinement has been considered, yielding a mesh that contains $240\times 880$ elements.\\
\begin{figure*}[!b]
\begin{center}
\begin{minipage}{0.45\textwidth}
\includegraphics[width=1.0\textwidth]{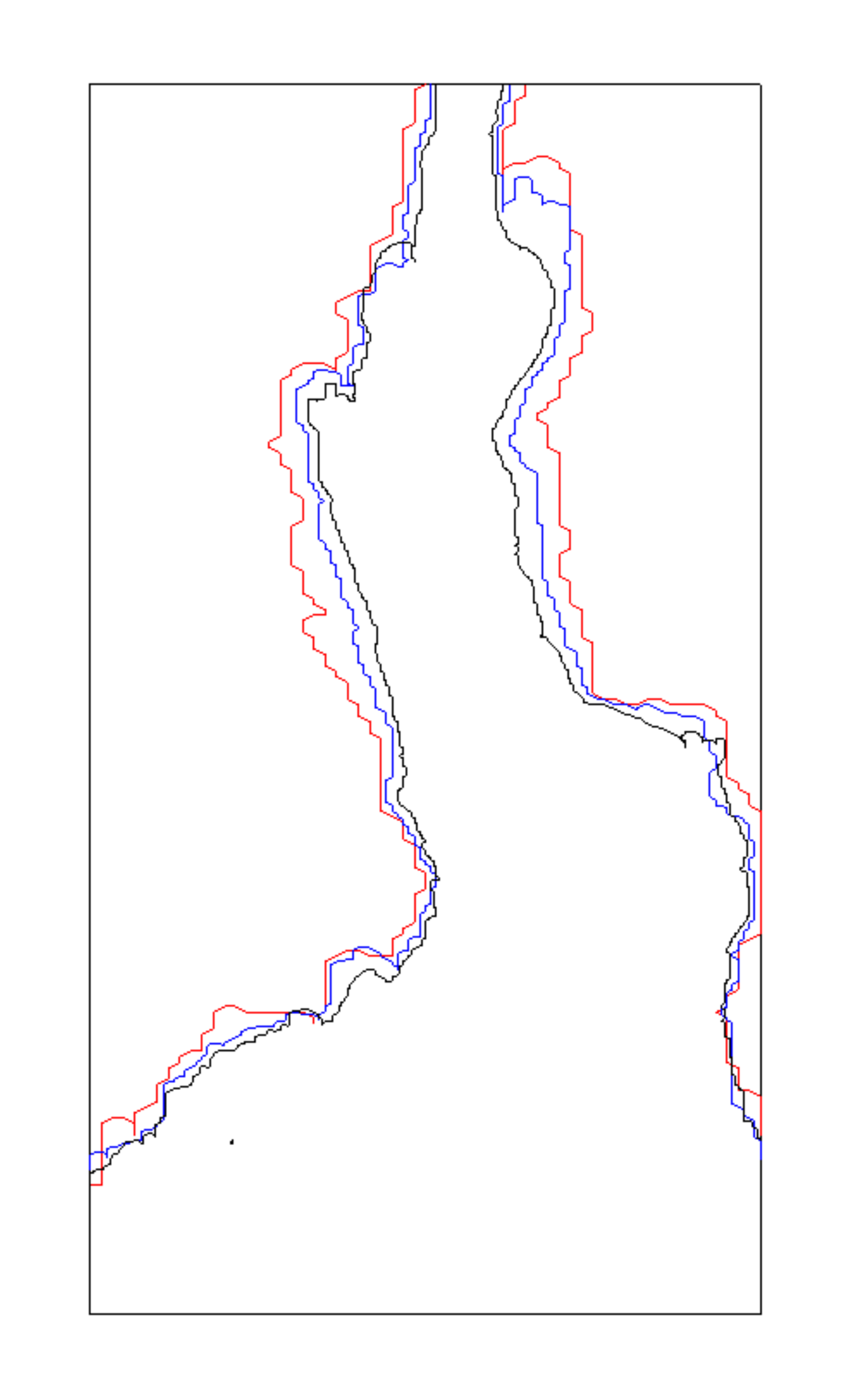}
\end{minipage}\hspace{-3mm}
\begin{minipage}{0.1\textwidth}
\includegraphics[width=1.0\textwidth]{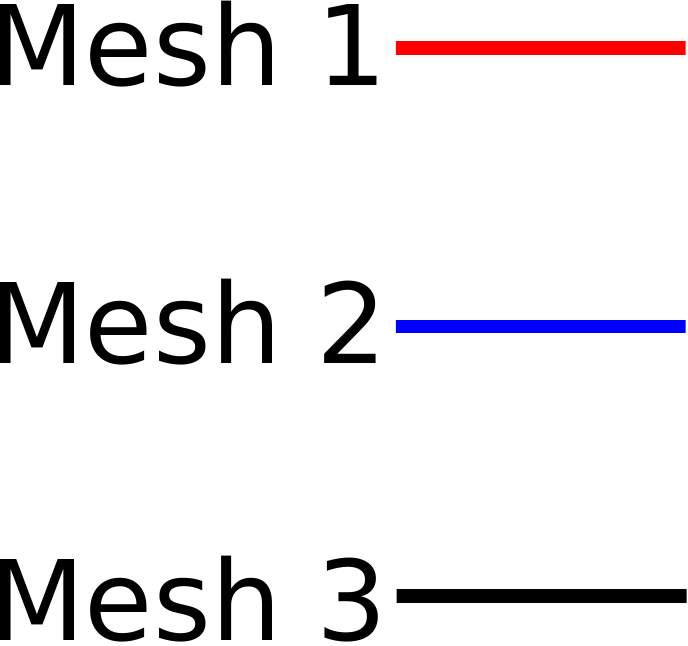}
\end{minipage}\hspace{-3mm}
\begin{minipage}{0.45\textwidth}
\includegraphics[width=1.0\textwidth]{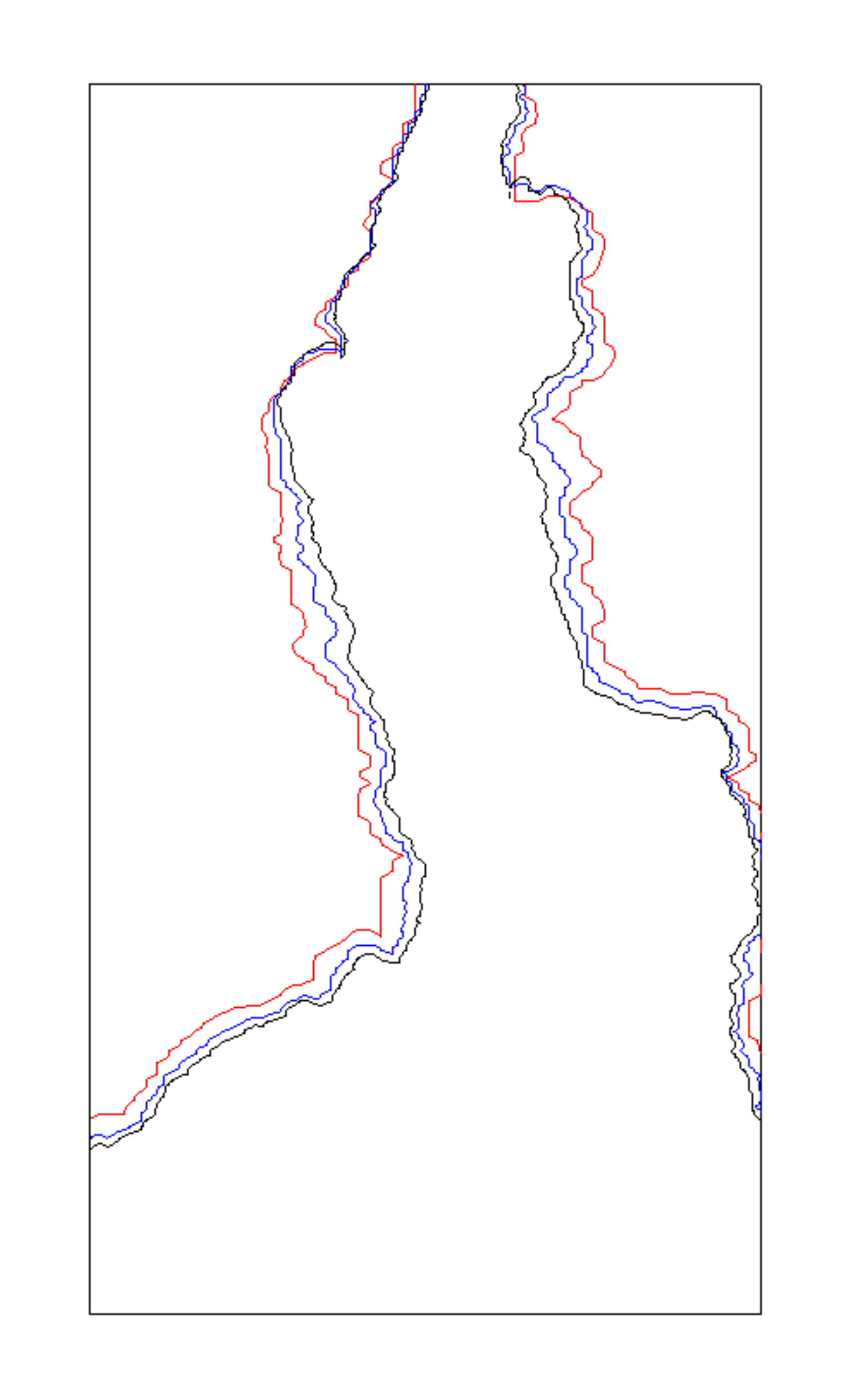}
\end{minipage}
\caption{Contour plots for $S_w=10^{-5}$ for the solutions at the end time $1.5\cdot 10^8\;[s]s$ obtained with the method developed in this paper (on the left) and for the fully-implicit TPFA solutions (on the right). In the legend, the corresponding mesh is provided \label{fig:contours2d}}
\end{center}
\end{figure*}
In Fig. \ref{fig:spe_solutions}, on one hand, we observe that all solutions possess the same behaviour in terms of front propagation and direction of flow. On the other hand, our method yield the same detection time for the water at the upper boundary independently of the mesh refinement, while a slower front propagation can be noticed for the first two fully-implicit TPFA solutions. The detection time for the first two fully-implicit TPFA solutions is $1.4\cdot 10^8\;[s]$. These differences are due to numerical diffusion, which cause a slower front propagation. In Fig. \ref{fig:contours2d}, contour plots for $S_w=10^{-5}$ at end time are presented for each method. Here, we can observe a reduced numerical diffusion of our method in comparison to the fully-implicit TPFA method. Therefore, the numerical diffusion causes the fully-implicit TPFA solution to yield a different detection time on the first two meshes, which can be properly reduced only on the finest mesh. 

\subsection{Anisotropic permeability}
\begin{figure*}[!hb]
\begin{center}
\begin{minipage}{0.32\textwidth}
\includegraphics[width=1.0\textwidth]{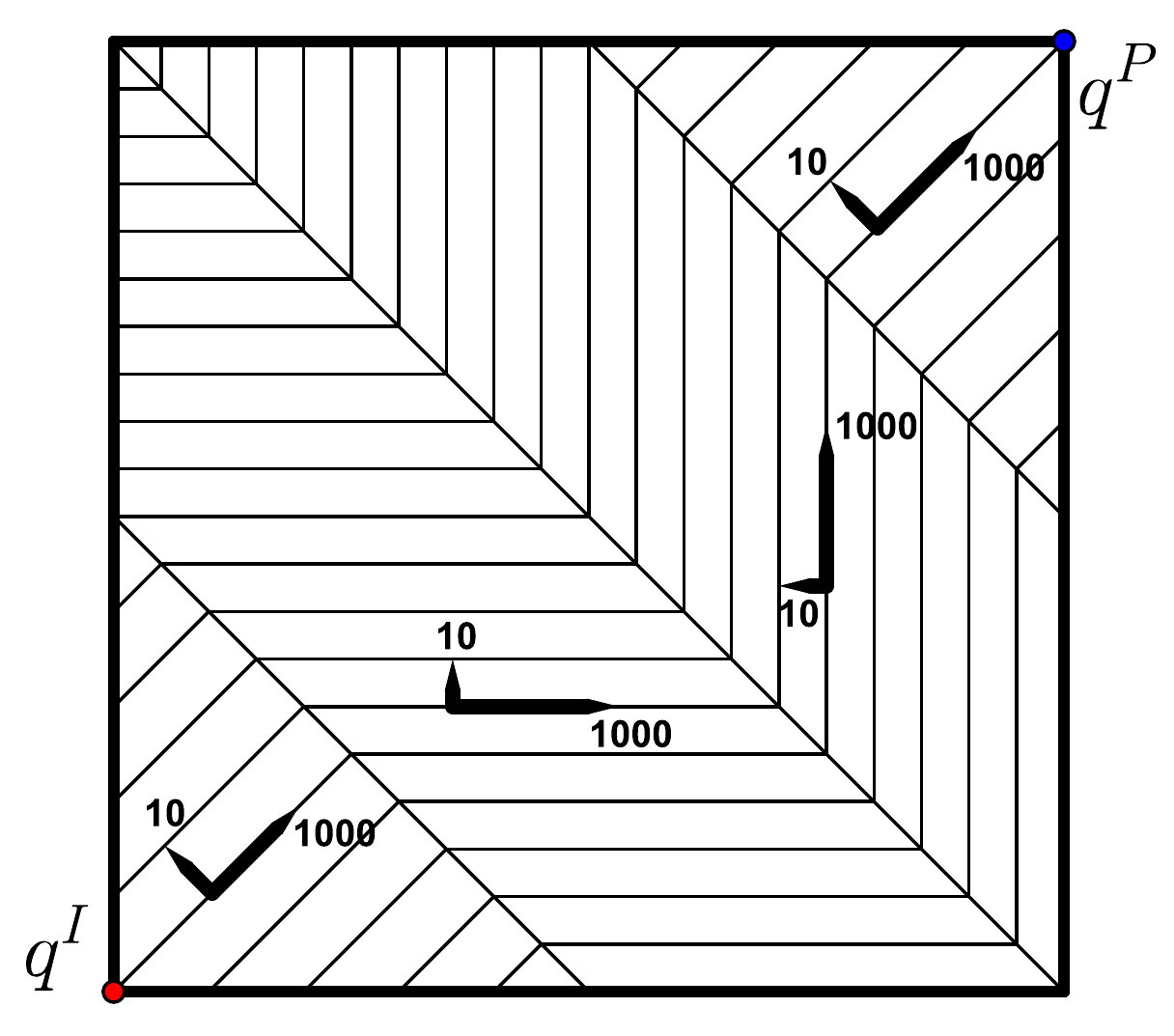}
\end{minipage}
\begin{minipage}{0.35\textwidth}
\includegraphics[width=1.0\textwidth]{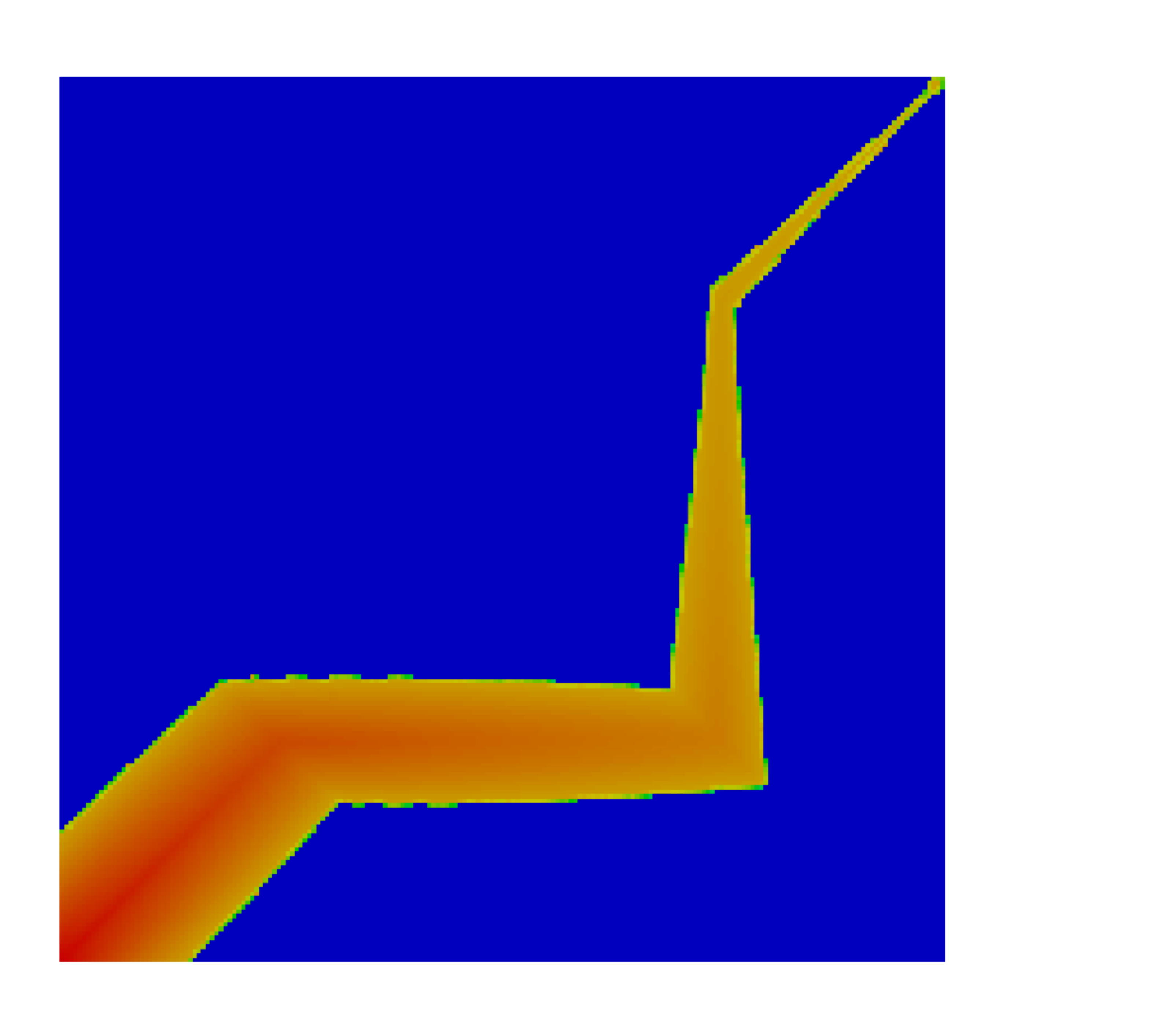}
\end{minipage}\hspace{-0.6cm}
\begin{minipage}{0.35\textwidth}
\includegraphics[width=1.0\textwidth]{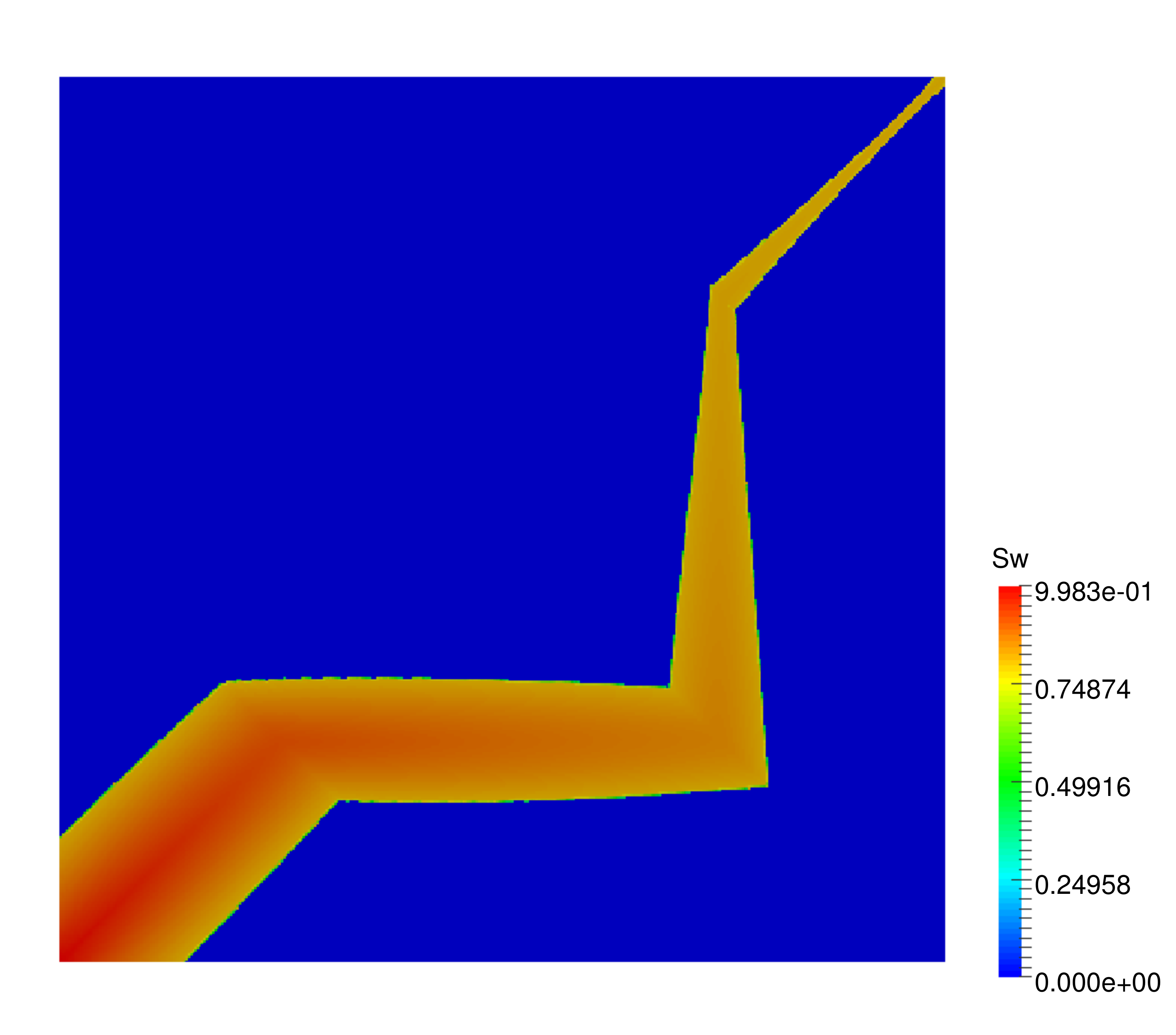}
\end{minipage}
\caption{On the left, the discontinuous anisotropic tensor and the positions of the injector and the producer are depicted. In the middle, the numerical solution on the coarse mesh is shown at end time, while on the right picture, the solution has been computed on a uniformly refined mesh \label{fig:anisotropic_perm}}
\end{center}
\end{figure*}

When incorporating realistic geological models, the numerical method has to be able to handle full-tensor permeabilities. This is investigated in the following example, which has been introduced in \cite{nikitin2014monotone}. Let $\fett{R}(\theta)$ denote the rotation matrix of angle $\theta$. Thus, the permeability (see Fig. \ref{fig:anisotropic_perm}, left) is defined as
$$
\fett{K}=10^{-13}\cdot\fett{R}(-\theta)\begin{pmatrix}
1000 & 0 \\ 0 & 10
\end{pmatrix}\fett{R}(\theta)\;[m^2],
$$
where the angle $\theta$ is equal to $45^\circ$ in the regions containing the wells, and alternates between $0^\circ$ and $90^\circ$ elsewhere. The problem domain is given as $\Omega=(0,100)^2$. The penalty parameter in \eqref{eq:penalty_parameter} is $\beta=10$. The viscosities are again $\mu_w=10^{-3}\;[Pa\cdot s]$ and $\mu_n=5.7\cdot 10^{-4}\;[Pa\cdot s]$. The flow is driven by the injection well $q^I$, located at the origin $(0,0)$, and the production well $q^P$, located at the upper-right corner $(100,100)$. The injection and production rates are given by
$$
\int_\Omega q^I = \int_\Omega q^P = 0.1\;[m^2/s],
$$
and no flow boundary condition is imposed. Gravity is neglected. The porosity is chosen constant on the entire domain, $\phi=0.2$. The relative permeabilities are chosen accordingly to \eqref{eq:BrooksCoreyLaw} with $\lambda=2$. The transport is simulated for $2.8\cdot 10^3$ seconds. The solution at the end time is depicted in Fig. \ref{fig:anisotropic_perm}, where the domain has been discretized by $200\times 200$ elements (middle) and by $400\times 400$ elements (right). We notice that the anisotropy is well captured and the amount of numerical diffusion is negligible. Furthermore, the front propagation of the wetting phase is the same for both simulations.

\subsection{Application to the three-dimensional SPE10 Benchmark}\label{sect:spe10-3d}
\begin{figure*}[!b]
\begin{center}
\includegraphics[width=0.32\textwidth]{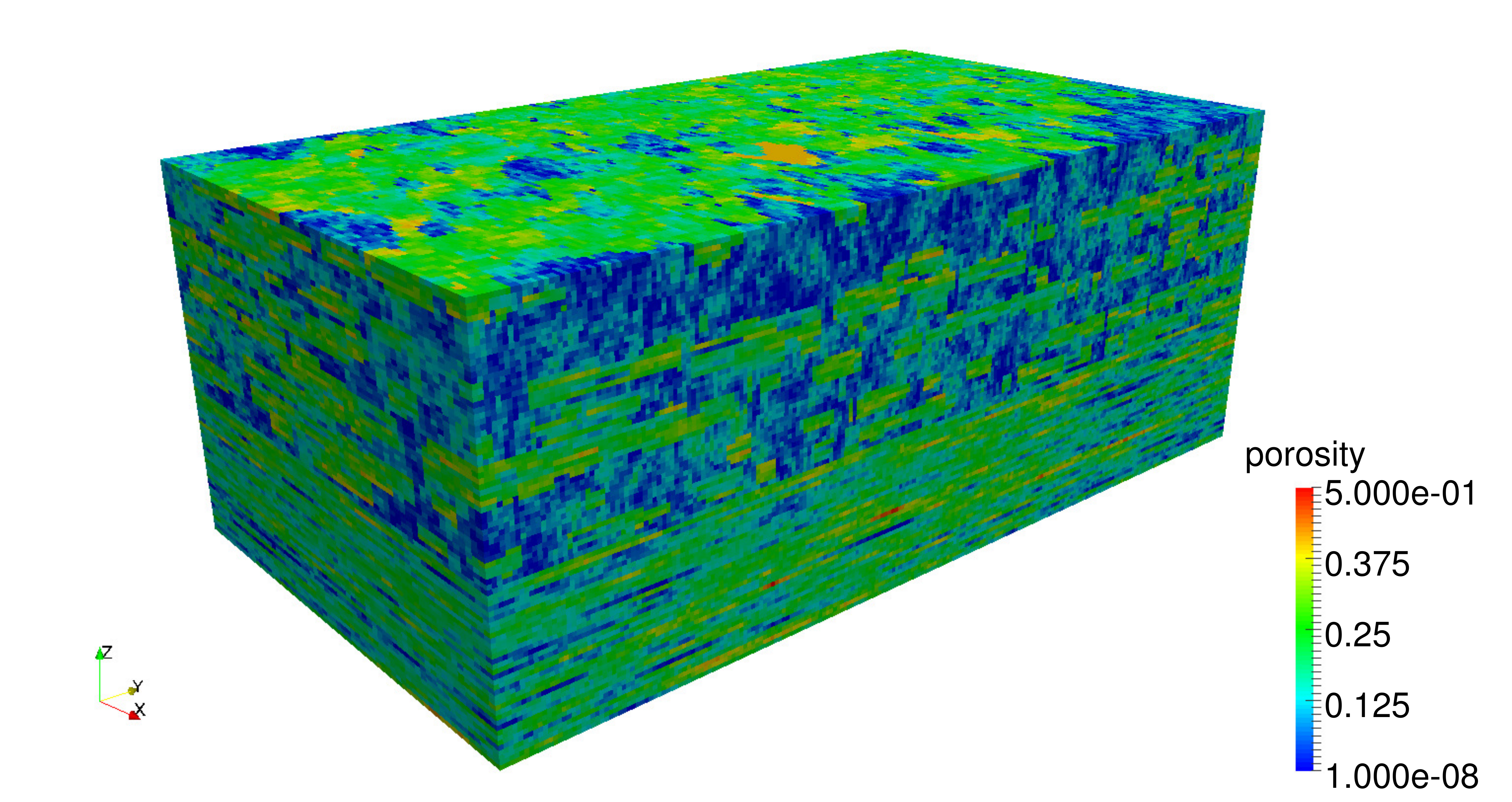}
\includegraphics[width=0.32\textwidth]{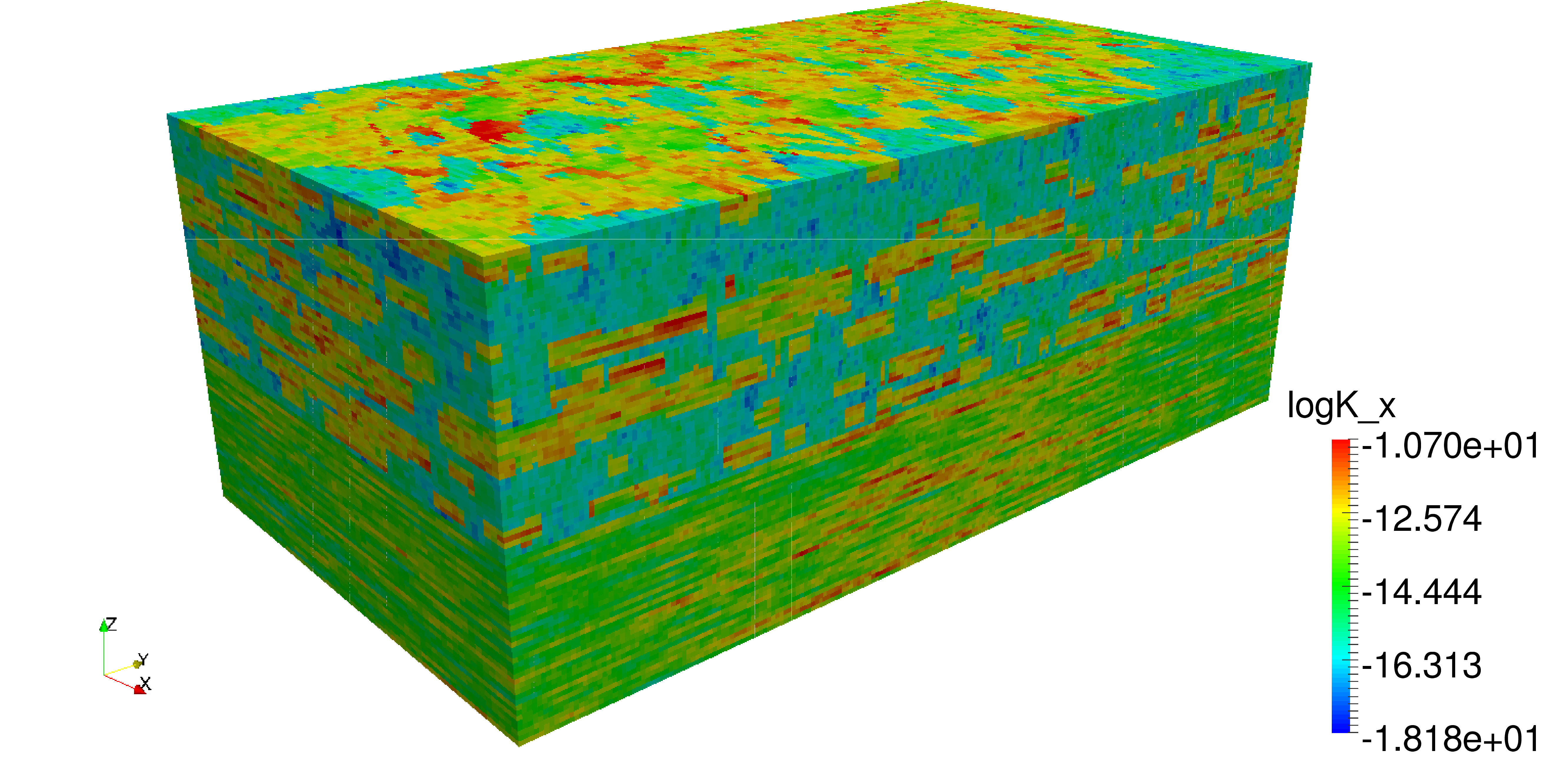}
\includegraphics[width=0.32\textwidth]{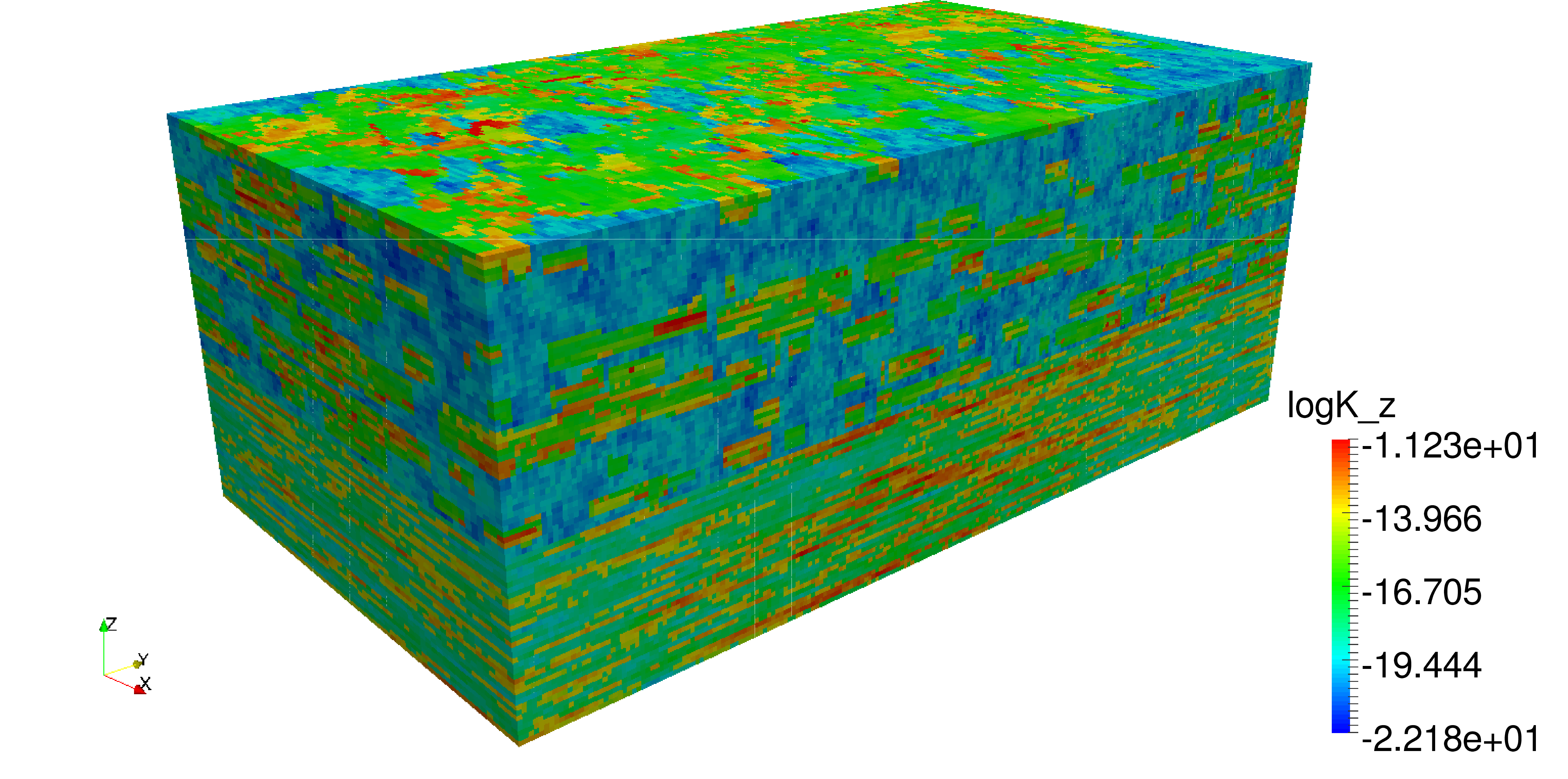}
\caption{On the left, the porosity is shown, while on the middle and on the right, the permeabilities in $x$ and $z$ direction are depicted, where the $z$-axis points in the depth direction. Permeabilities are isotropic in the $x$ and $y$ direction, i.e. $K_x = K_y$ ($z$-direction scaled by 5) \label{fig:porosity_perm_spe10_3d}}
\end{center}
\end{figure*}

\begin{figure*}[!hb]
\begin{center}
\begin{minipage}{0.5\textwidth}
\includegraphics[width=1.0\textwidth]{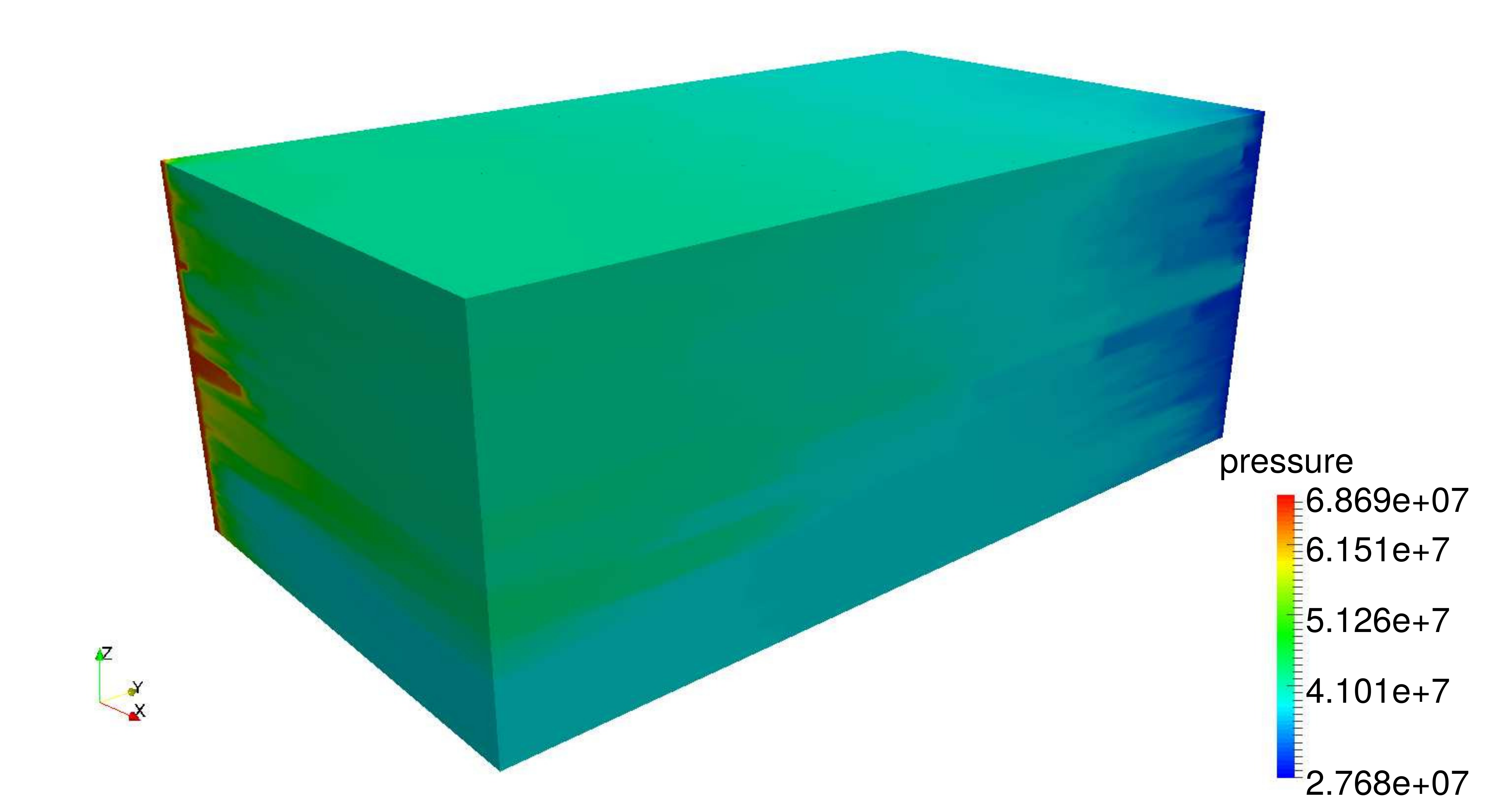}
\end{minipage}\hspace{-0.3cm}
\begin{minipage}{0.5\textwidth}
\includegraphics[width=1.0\textwidth]{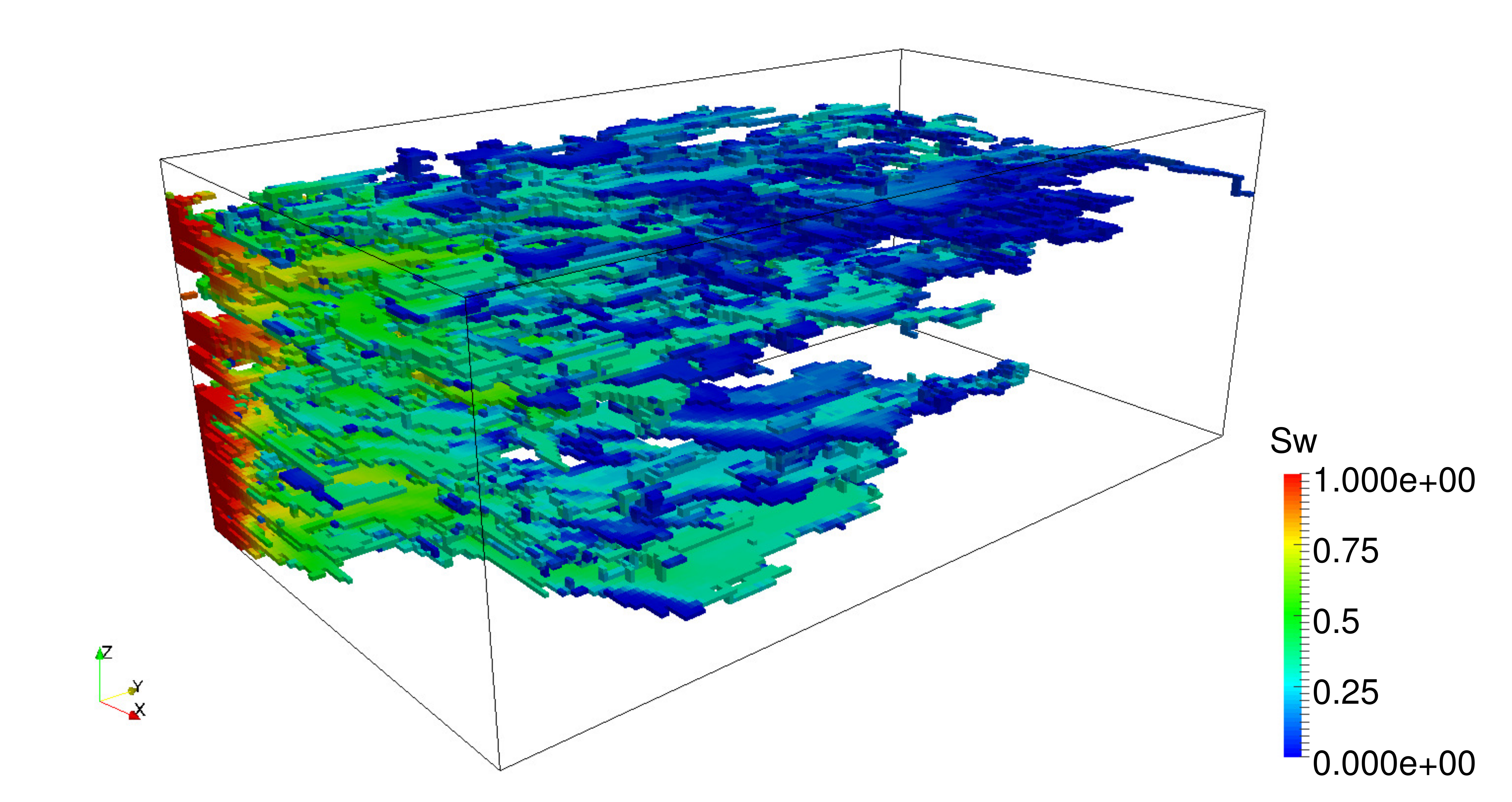}
\end{minipage}
\caption{On the left, the numerical solution of the pressure field for the SPE10 Benchmark is shown. On the right, saturation for the wetting phase is depicted over a threshold of 0.001 after 240 days ($z$-direction scaled by 5) \label{fig:spe10_3d_sol}}
\end{center}
\end{figure*}

In the previous numerical tests, we tested the method developed in this paper for two-dimensional problems, where it was demonstrated that fronts are well captured by our method with less numerical diffusion compared to the standard fully-implicit TPFA method. In addition, it was shown that the scheme is also consistent for full anisotropic tensors. In the following example, we apply our method on a well-established three-dimensional problem with realistic geological data. The setting of this example is taken from the second problem of the SPE10 Benchmark \cite{christie2001tenth}. The domain $\Omega$ is discretized by $60\times 220\times 85$ cells and has a size of $L_x\times L_y\times L_z = 365.76\times 670.56\times 51.816\;[m]$. This yields the discretization sizes $dx=6.096\;[m]$, $dy=3.048\;[m]$ and $dz=0.6096\;[m]$. The permeability and porosity fields are depicted in Fig. \ref{fig:porosity_perm_spe10_3d}. The model consists of two different geological formations: a shallow-marine Tarbert formation in the top 35 layers and a fluivial Upper-Ness formation in the bottom 50 layers. In the upper formation, the permeability is relatively smooth, while the bottom formation possesses a more heterogeneous structure, including channels. In both formations, the permeabilities are characterized by large variations of 8-12 orders of magnitude.\\
We are interested in validating our method in terms of front propagation, and therefore we consider a simplified version of the original SPE10 Benchmark problem.  Here, we simulate a flow driven by a pressure gradient, i.e., we impose Dirichlet boundary condition on the following sets:
$$
\Gamma_{D,\mathrm{inj}}= \left\{ (x,y,z)\in \partial\Omega : x<dx, y<dy \right\},
$$
$$
\Gamma_{D,\mathrm{ext}}= \left\{ (x,y,z)\in \partial\Omega : x>L_x-dx, y>L_y-dy \right\}.
$$
On $\Gamma_{D,\mathrm{inj}}$, we set $g_D=6.8948\cdot 10^7\;[Pa]$, while on $\Gamma_{D,\mathrm{ext}}$, we impose $g_D=2.7579\cdot 10^7\;[Pa]$. On the rest of the boundary $\partial\Omega\setminus \left( \Gamma_{D,\mathrm{inj}} \cup \Gamma_{D,\mathrm{ext}}\right)$, no-flow condition is provided. Furthermore, we assume the wetting phase (water) infiltrates from $\Gamma_{D,\mathrm{inj}}$, i.e., $S_w=1$.\\
Viscosity for the wetting phase is $\mu_w = 3\cdot 10^{-4}\;[Pa\cdot s]$, and for the non-wetting phase (oil) is $\mu_n = 3\cdot 10^{-3}\;[Pa\cdot s]$.
The relative permeabilities are chosen accordingly to the quadratic law \eqref{eq:quadratic_rel_perm}. Initially, the domain is filled by the non-wetting phase. We choose a uniform time step of 10 days and the simulation runs until the water reaches the boundary $\Gamma_{D,\mathrm{ext}}$.\\
The original mesh consists of 1,220,000 elements making it hard to solve sequentially the linear system obtained with the SWIPG method. Therefore, the solution has been computed on 35 processors, yielding a total number of 1,326,960 elements, which includes the overlapping elements. For visualization purposes the $z$-axis is scaled by a factor five in all figures showing the three-dimensional model domain.\\
We obtained a detection time of 240 days for the water front. The wetting phase saturation at this time is shown in Fig. \ref{fig:spe10_3d_sol}, on the right.
On the left side of Fig. \ref{fig:spe10_3d_sol}, the pressure field after the first time step is provided. As reference, a numerical solution is computed using a fully-implicit  TPFA method. The reference solution yields a detection time of 250 days, which is in good agreement with the result obtained with our method (deviation of $\Delta t$). As in the two-dimensional heterogeneous problem from Sect. \ref{sect:2d_heterogeneous_prob}, the difference in the front propagation is caused by a higher numerical diffusion of the fully-implicit TPFA method.
\section{Final remarks}\label{sect:conclusion}
In this work, we have proposed an improved streamline approach for the fast simulation of incompressible two-phase flow in porous media for high-rate flooding scenarios, where capillary effects can be neglected. We have applied a sequential algorithm, where the pressure equation is solved by a DG method, while the system of one-dimensional Riemann problems along streamlines is solved using the front tracking method. The advantages of our method consist in combining the optimal approximation properties of the DG method with a fast and unconditionally stable solver for the transport equation. Furthermore, a parallel version of the algorithm for the streamline tracing on decomposed domains has been presented. A series of numerical tests for two- and three-dimensional problems has shown the reliability of the presented method in terms of flow front approximation. The reduced diffusivity of our method has also been shown in comparison to a standard fully-implicit TPFA method for different two-dimensional examples. Future work will include further improvements of the method in the direction of a better mass-conservation. An extension to non-cartesian grids will also be investigated, in order to include more complicated geometries.

\section*{Acknowledgements}
This work was partially supported by the DFG grant (WO/671 11-1)
\bibliographystyle{siam}
\bibliography{vidotto_preprint}

\begin{thebibliography}{10}

\bibitem{aavatsmark2007interpretation}
{\sc I.~Aavatsmark}, {\em Interpretation of a two-point flux stencil for skew
  parallelogram grids}, Computational geosciences, 11 (2007), p.~199.

\bibitem{ask2000local}
{\sc A.~Ask, H.~K. Dahle, K.~H. Karlsen, and H.~F. Nordhaug}, {\em A local
  streamline eulerian-lagrangian method for two-phase flow}, The XIII
  international conference on computational methods in water resources, 2000.

\bibitem{bastian2011benchmark}
{\sc P.~Bastian}, {\em Benchmark {3D}: Symmetric weighted interior penalty
  discontinuous {Galerkin} scheme}, in Finite Volumes for Complex Applications
  VI Problems \& Perspectives, Springer, 2011, pp.~949--959.

\bibitem{bastian2014fully}
\leavevmode\vrule height 2pt depth -1.6pt width 23pt, {\em A fully-coupled
  discontinuous {Galerkin} method for two-phase flow in porous media with
  discontinuous capillary pressure}, Computational Geosciences, 18 (2014),
  pp.~779--796.

\bibitem{bastian2008generic}
{\sc P.~Bastian, M.~Blatt, A.~Dedner, C.~Engwer, R.~Kl{\"o}fkorn, R.~Kornhuber,
  M.~Ohlberger, and O.~Sander}, {\em A generic grid interface for parallel and
  adaptive scientific computing. {Part II}: Implementation and tests in
  {DUNE}}, Computing, 82 (2008), pp.~121--138.

\bibitem{bastian2012algebraic}
{\sc P.~Bastian, M.~Blatt, and R.~Scheichl}, {\em Algebraic multigrid for
  discontinuous {Galerkin} discretizations of heterogeneous elliptic problems},
  Numerical Linear Algebra with Applications, 19 (2012), pp.~367--388.

\bibitem{bastian2003superconvergence}
{\sc P.~Bastian and B.~Rivi{\`e}re}, {\em Superconvergence and {H(div)}
  projection for discontinuous {Galerkin} methods}, International journal for
  numerical methods in fluids, 42 (2003), pp.~1043--1057.

\bibitem{batycky19973d}
{\sc R.~Batycky, M.~J. Blunt, M.~R. Thiele, et~al.}, {\em A {3D} field-scale
  streamline-based reservoir simulator}, SPE Reservoir Engineering, 12 (1997),
  pp.~246--254.

\bibitem{berre2002streamline}
{\sc I.~Berre, H.~K. Dahle, K.~H. Karlsen, and H.~F. Nordhaug}, {\em A
  streamline front tracking method for two-and three-phase flow including
  capillary forces}, Contemporary Mathematics, 295 (2002), pp.~49--62.

\bibitem{bhambri2011compositional}
{\sc P.~Bhambri and K.~Mohanty}, {\em Compositional streamline simulation: a
  parallel implementation}, Transport in porous media, 90 (2011), pp.~741--761.

\bibitem{borah2013investigation}
{\sc A.~Borah, P.~Singh, and P.~Goswami}, {\em An {Investigation} of {Solving
  Multidimensional Multiphase Flow}: {Streamline} front tracking method},
  (2013).

\bibitem{bratvedt1996streamline}
{\sc F.~Bratvedt, T.~Gimse, and C.~Tegnander}, {\em Streamline computations for
  porous media flow including gravity}, Transport in Porous Media, 25 (1996),
  pp.~63--78.

\bibitem{brezzi1985two}
{\sc F.~Brezzi, J.~Douglas, and L.~D. Marini}, {\em Two families of mixed
  finite elements for second order elliptic problems}, Numerische Mathematik,
  47 (1985), pp.~217--235.

\bibitem{brezzi2012mixed}
{\sc F.~Brezzi and M.~Fortin}, {\em Mixed and hybrid finite element methods},
  vol.~15, Springer Science \& Business Media, 2012.

\bibitem{camp2011streamline}
{\sc D.~Camp, C.~Garth, H.~Childs, D.~Pugmire, and K.~Joy}, {\em Streamline
  integration using {MPI}-hybrid parallelism on a large multicore
  architecture}, IEEE Transactions on Visualization and Computer Graphics, 17
  (2011), pp.~1702--1713.

\bibitem{cao2010robust}
{\sc Y.~Cao}, {\em Robust {Numerical} {Algorithms} {Based} on {Corrected}
  {Operator} {Splitting} for {Two}-phase {Flow} in {Porous} {Media}}, Shaker
  Verlag Gmbh, 2010.

\bibitem{cao2007fractional}
{\sc Y.~Cao, B.~Eikemo, and R.~Helmig}, {\em Fractional flow formulation for
  two-phase flow in porous Media}, GRK 1398/1, 2007.

\bibitem{cao2011two}
{\sc Y.~Cao, R.~Helmig, and B.~Wohlmuth}, {\em A two-scale operator-splitting
  method for two-phase flow in porous media}, Advances in Water Resources, 34
  (2011), pp.~1581--1596.

\bibitem{christie2001tenth}
{\sc M.~Christie, M.~Blunt, et~al.}, {\em Tenth {SPE} comparative solution
  project: {A} comparison of upscaling techniques}, in SPE Reservoir Simulation
  Symposium, Society of Petroleum Engineers, 2001.

\bibitem{crane1999streamline}
{\sc M.~J. Crane and M.~J. Blunt}, {\em Streamline-based simulation of solute
  transport}, Water Resources Research, 35 (1999), pp.~3061--3078.

\bibitem{dake2001practice}
{\sc L.~P. Dake}, {\em The practice of reservoir engineering (revised
  edition)}, vol.~36, Elsevier, 2001.

\bibitem{datta2007streamline}
{\sc A.~Datta-Gupta and M.~J. King}, {\em Streamline simulation: {Theory} and
  practice}, vol.~11, Society of Petroleum Engineers Richardson, 2007.

\bibitem{di2012analysis}
{\sc D.~A. Di~Pietro and A.~Ern}, {\em Analysis of a discontinuous galerkin
  method for heterogeneous diffusion problems with low-regularity solutions},
  Numerical Methods for Partial Differential Equations, 28 (2012),
  pp.~1161--1177.

\bibitem{epshteyn2008convergence}
{\sc Y.~Epshteyn and B.~Riviere}, {\em Convergence of high order methods for
  miscible displacement}, International Journal of Numerical Analysis and
  Modeling, 5 (2008), pp.~47--63.

\bibitem{gerritsen2009parallel}
{\sc M.~G. Gerritsen, H.~L{\"o}f, and M.~R. Thiele}, {\em Parallel
  implementations of streamline simulators}, Computational Geosciences, 13
  (2009), p.~135.

\bibitem{helmig1997multiphase}
{\sc R.~Helmig}, {\em Multiphase flow and transport processes in the
  subsurface: a contribution to the modeling of hydrosystems.},
  Springer-Verlag, 1997.

\bibitem{holden2010splitting}
{\sc H.~Holden}, {\em Splitting methods for partial differential equations with
  rough solutions: {Analysis} and {MATLAB} programs}, European Mathematical
  Society, 2010.

\bibitem{holden2015front}
{\sc H.~Holden and N.~H. Risebro}, {\em Front tracking for hyperbolic
  conservation laws}, vol.~152, Springer, 2015.

\bibitem{kane2017hp}
{\sc B.~Kane, R.~Kl{\"o}fkorn, and C.~Gersbacher}, {\em hp--{Adaptive}
  {Discontinuous} {Galerkin} {Methods} for {Porous} {Media} {Flow}}, in
  International Conference on Finite Volumes for Complex Applications,
  Springer, 2017, pp.~447--456.

\bibitem{kippe2007method}
{\sc V.~Kippe, H.~H{\ae}gland, K.-A. Lie, et~al.}, {\em A method to improve the
  mass balance in streamline methods}, in SPE Reservoir Simulation Symposium,
  Society of Petroleum Engineers, 2007.

\bibitem{langseth1996implementation}
{\sc J.~Langseth}, {\em On an implementation of a front tracking method for
  hyperbolic conservation laws}, Advances in engineering software, 26 (1996),
  pp.~45--63.

\bibitem{li2015high}
{\sc J.~Li and B.~Riviere}, {\em High order discontinuous {Galerkin} method for
  simulating miscible flooding in porous media}, Computational Geosciences, 19
  (2015), p.~1251.

\bibitem{li2004tidal}
{\sc L.~Li, D.~Barry, D.~Jeng, and H.~Prommer}, {\em Tidal dynamics of
  groundwater flow and contaminant transport in coastal aquifers}, Coastal
  aquifer management: Monitoring, modeling, and case studies,  (2004),
  pp.~115--141.

\bibitem{lin2015comparative}
{\sc G.~Lin, J.~Liu, and F.~Sadre-Marandi}, {\em A comparative study on the
  weak {Galerkin}, discontinuous {Galerkin}, and mixed finite element methods},
  Journal of Computational and Applied Mathematics, 273 (2015), pp.~346--362.

\bibitem{mallison2006improved}
{\sc B.~T. Mallison, M.~G. Gerritsen, S.~F. Matringe, et~al.}, {\em Improved
  mappings for streamline-based simulation}, SPE Journal, 11 (2006),
  pp.~294--302.

\bibitem{niessner2007multi}
{\sc J.~Niessner and R.~Helmig}, {\em Multi-scale modeling of
  three-phase--three-component processes in heterogeneous porous media},
  Advances in Water Resources, 30 (2007), pp.~2309--2325.

\bibitem{nikitin2014monotone}
{\sc K.~Nikitin, K.~Terekhov, and Y.~Vassilevski}, {\em A monotone nonlinear
  finite volume method for diffusion equations and multiphase flows},
  Computational Geosciences, 18 (2014), pp.~311--324.

\bibitem{nilsen2009front}
{\sc H.~M. Nilsen, K.-A. Lie, et~al.}, {\em Front {Tracking} {Methods} for
  {Use} in {Streamline} {Simulation} of {Compressible} {Flow}}, in SPE
  Reservoir Simulation Symposium, Society of Petroleum Engineers, 2009.

\bibitem{oden1998discontinuoushpfinite}
{\sc J.~T. Oden, I.~Babu{\^s}ka, and C.~E. Baumann}, {\em A discontinuous hp
  finite element method for diffusion problems}, Journal of computational
  physics, 146 (1998), pp.~491--519.

\bibitem{petruzzelli2013migration}
{\sc D.~Petruzzelli and F.~G. Helfferich}, {\em Migration and fate of
  pollutants in soils and subsoils}, vol.~32, Springer Science \& Business
  Media, 2013.

\bibitem{pollock1988semianalytical}
{\sc D.~W. Pollock}, {\em Semianalytical computation of path lines for
  finite-difference models}, Ground water, 26 (1988), pp.~743--750.

\bibitem{pugmire2009scalable}
{\sc D.~Pugmire, H.~Childs, C.~Garth, S.~Ahern, and G.~H. Weber}, {\em Scalable
  computation of streamlines on very large datasets}, in Proceedings of the
  Conference on High Performance Computing Networking, Storage and Analysis,
  ACM, 2009, p.~16.

\bibitem{riviere2008discontinuous}
{\sc B.~Riviere}, {\em Discontinuous {Galerkin} methods for solving elliptic
  and parabolic equations: theory and implementation}, Society for Industrial
  and Applied Mathematics, 2008.

\bibitem{schneider2016monotone}
{\sc M.~Schneider, B.~Flemisch, and R.~Helmig}, {\em Monotone nonlinear
  finite-volume method for nonisothermal two-phase two-component flow in porous
  media}, International Journal for Numerical Methods in Fluids,  (2016).

\bibitem{siavashi2014three}
{\sc M.~Siavashi, M.~J. Blunt, M.~Raisee, and P.~Pourafshary}, {\em
  Three-dimensional streamline-based simulation of non-isothermal two-phase
  flow in heterogeneous porous media}, Computers \& Fluids, 103 (2014),
  pp.~116--131.

\bibitem{vasco1998integrating}
{\sc D.~W. Vasco, S.~Yoon, A.~Datta-Gupta, et~al.}, {\em Integrating dynamic
  data into high-resolution reservoir models using streamline-based analytic
  sensitivity coefficients}, in SPE Annual Technical Conference and Exhibition,
  Society of Petroleum Engineers, 1998.

\bibitem{zunino2009discontinuous}
{\sc P.~Zunino}, {\em Discontinuous {Galerkin} methods based on weighted
  interior penalties for second order {PDEs} with non-smooth coefficients},
  Journal of Scientific Computing, 38 (2009), pp.~99--126.

\end{thebibliography}

\end{document}